# Improved Optimization of Finite Sums with Minibatch Stochastic Variance Reduced Proximal Iterations

Jialei Wang[*] and Tong Zhang[♯]

[*]Department of Computer Science, University of Chicago, IL, USA
[♯]Tencent AI Lab, Shenzhen, China

October 10, 2017

**Abstract**

We present a novel minibatch stochastic optimization method for empirical risk minimization problems. The method efficiently leverages variance reduced first-order and sub-sampled higher-order information to accelerate the convergence. We prove improved iteration complexity over state-of-the-art methods under suitable conditions. Experiments are provided to demonstrate the empirical advantage of the proposed method over existing approaches in the literature.

## 1 Introduction

We consider the following optimization problem of finite-sums:

$$\min_w f(w) := \frac{1}{n}\sum_{i=1}^{n} f_i(w), \tag{1}$$

which arises in many machine learning problems. In the context of regularized loss minimization of linear predictors, we have $f_i(w) = \ell(w^\top x_i; b_i) + \frac{\lambda}{2}\|w\|^2$. Let $x_1, x_2, ..., x_n \in \mathbb{R}^d$ be feature vectors of $n$ data samples, and $b_1, ..., b_n \in \mathbb{R}$ or $\{-1, 1\}$ be the corresponding target variables of interest. We note that (1) covers many popular models used in machine learning: for example when $\ell(w^\top x_i; b_i) = \frac{1}{2}(b_i - \langle w, x_i\rangle)^2$ we obtain ridge regression; when $\ell(w^\top x_i; b_i) = \log(1+\exp(-b_i\langle w, x_i\rangle))$ we obtain $\ell_2$ regularized logistic regression. Other examples such as SVMs can also be obtained by setting different loss functions $\ell(w^\top x_i; b_i)$. The ubiquitousness of such finite-sum optimization problems and the massive scale of modern datasets motivate significant research effort on efficient optimization algorithms to solve (1).

In this paper, we let

$$w^* = \arg\min_w f(w)$$

be the optimum of (1). For any approximate solution $w$, we say that it achieves $\epsilon$-objective suboptimality if $f(w) - f(w^*) \leq \epsilon$. For any $\epsilon > 0$, the runtime cost of optimization algorithms to find an $\epsilon$-suboptimal minimizer typically depends on the target accuracy $\epsilon$, as well as the conditions of the problem (1). Throughout the paper, we use the following notion of strong convexity and smoothness to discuss the iteration complexity of various optimization algorithms.

**Definition 1.** *A function $f_i(w)$ is $L$-smooth with respect to $w$ if $f_i(w)$ is differentiable and its gradient is $L$-Lipschitz continuous, i.e. we have*

$$\left\|\nabla f_i(w) - \nabla f_i(w')\right\| \le L \left\|w - w'\right\|, \quad \forall w, w' \in \mathbb{R}^d.$$

*A function $f_i(w)$ is $\lambda$-strongly convex with respect to $w$, if we have*

$$f(w) \ge f(w') + \langle \nabla f(w'), w - w' \rangle + \frac{\lambda}{2} \left\|w - w'\right\|^2, \quad \forall w, w' \in \mathbb{R}^d.$$

It is well-known that a consequence of $L$-smoothness is the following quadratic upper bound for $f(w)$:

$$f(w) \le f(w') + \langle \nabla f(w'), w - w' \rangle + \frac{L}{2} \left\|w - w'\right\|^2, \quad \forall w, w' \in \mathbb{R}^d.$$

In recent years there have been significant advances in developing fast optimization algorithms for (1), and we refer the readers to (Bottou et al., 2016) for a comprehensive survey of these developments. For large-scale problems in the form of (1), randomized methods are particularly efficient because of their low per iteration computation. Below we briefly review two lines of research: i) randomized variance reduced first-order methods; ii) randomized methods leveraging second-order information. Through out this section we focus on problems where each function $f_i(w), \forall i \in [n]$ is $\lambda$-strongly convex and $L$-smooth. For non-strongly convex or non-smooth objectives, we can use regularization or smoothing technique to reduce it to strongly convex and smooth objectives (Allen-Zhu and Hazan, 2016).

## 1.1 Variance Reduced First-order Methods

The key technique for developing fast stochastic first-order methods is variance reduction, which makes the variance of the randomized update direction approching zero when the iterate gets closer to the optimum. Representative methods of this category include SAG (Roux et al., 2012), SVRG (Johnson and Zhang, 2013), SDCA (Shalev-Shwartz and Zhang, 2013) and SAGA (Defazio et al., 2014), etc. In order to find a solution that reaches $\epsilon$-suboptimality, these methods requires

$$\mathcal{O}\left(\left(n + \frac{L}{\lambda}\right) \log\left(\frac{1}{\epsilon}\right)\right) \tag{2}$$

calls to the first-order oracle to compute the gradient of an individual function. In the large condition number regime (e.g. $\frac{L}{\lambda} > n$), by using an acceleration technique (such as Catalyst (Shalev-Shwartz and Zhang, 2016; Frostig et al., 2015; Lin et al., 2015a), APCG (Lin et al., 2015b), SPDC (Zhang and Xiao, 2015), Katyusha (Allen-Zhu, 2016), etc), one can further improve the iteration complexity to

$$\mathcal{O}\left(\left(n + \sqrt{n \cdot \frac{L}{\lambda}}\right) \log\left(\frac{1}{\epsilon}\right)\right). \tag{3}$$

Recently, Woodworth and Srebro (2016) provided a first-order oracle lower bound for (1) which has the same iteration complexity as (3). This means that (3) is not improvable for general (high dimensional) problems, if the optimization procedure only uses the first order and proximal oracles of individual functions.

### 1.2 Leveraging Second-order Information

Second-order information is often useful in improving the convergence of optimization algorithms. However for large-scale problems, obtaining and inverting the exact Hessian matrix is often computational expensive, making vanilla Newton methods not well suited in solving (1). Therefore, there have been emerging studies in designing randomized algorithms that effectively utilize the *approximated* Hessian information. One line of such research is the sub-sampled Newton method which approximates the Hessian matrix based on a sub-sampled minibatch of data. Several work, such as (Byrd et al., 2011; Erdogdu and Montanari, 2015; Roosta-Khorasani and Mahoney, 2016a,b; Bollapragada et al., 2016), established local linear convergence rates for different variants of sub-sampled Newton methods, when the size of the sub-sampled data is large enough. (Xu et al., 2016) considered non-uniform sampling in constructing the stochastic Hessian matrix, and showed that weighted sampling according to individual smoothness or leverage score can help constructing better approximation of Hessian, which improves the convergence. (Pilanci and Wainwright, 2017) discussed how to use sketching instead of sub-sampling to approximate the Hessian matrix. All of the aforementioned methods employ the full first-order gradient, and require calling the second-order oracle to compute the Hessian and solving a resulting linear system at every iteration. Therefore the computational complexity (as compared in Table 2 of (Xu et al., 2016)) are often worse than stochastic algorithms such as SVRG.

Another line of research is to consider randomness in both first and second-order information to design lower per-iteration cost algorithms. In particular, (Byrd et al., 2016) considered combining minibatch SGD with Limited-memory BFGS (L-BFGS) (Liu and Nocedal, 1989) type update. Inspired by this, (Moritz et al., 2016; Gower et al., 2016; Wang et al., 2017b) proposed to combine variance reduced stochastic gradient with L-BFGS, and proved linear convergence for strongly convex and smooth objectives. Moreover, experiments in (Moritz et al., 2016; Gower et al., 2016; Wang et al., 2017b) also demonstrated superior empirical performance of this type of methods. However, theoretically it is hard to guarantee the quality of approximated Hessian using L-BFGS update, and thus the iteration complexity obtained by (Moritz et al., 2016; Gower et al., 2016; Wang et al., 2017b) can be pessimistic, and can be much worse than vanilla SVRG. Moreover, (Qu et al., 2015) proposed to incorporate curvature information in minibatch SDCA methods, and showed improved convergence over SDCA, but the method proposed in (Qu et al., 2015) involves solving a much more expensive subproblem than minibatch SDCA. Therefore the overall runtime benefit is still unclear. (Gonen et al., 2016) considered ridge regression problems specifically, and proposed to use sketching to compute the rank-$k$ approximation of the Hessian matrix, based on which SVRG is ran with preconditioning. (Yang et al., 2016) suggested to use preconditioning as a preprocessing step for general stochastic optimization problems. (Agarwal et al., 2016) proposed to approximate the inverse Hessian matrix directly by sampling from its Taylor expansion. In terms of lower bound, (Arjevani and Shamir, 2016) showed under some mild algorithmic assumptions, second-order oracle generally cannot improve the oracle complexity over first-order methods, for very high-dimensional problems.

Though lower bound have been established showing that the iteration complexity of accelerated SVRG methods cannot be improved in general, the establishments of these lower bounds are based on constructing very high-dimensional hard problems where dimension $d$ can be much larger than sample size $n$. Such pessimistic situations are not consistent with practice, where second-order information is observed to be very helpful in improving the convergence speed. Therefore it is still interesting to analyze theoretically how second-order information can be helpful under suitable conditions, and more importantly, how to design more efficient methods that make use of possibly noisy second-order information. In this paper, we make effort in this direction and make the following contributions:

- We propose a novel approach that combines the advantages of variance-reduced first-order methods and sub-sampled Newton methods, in a efficient way that does not require expensive Hessian matrix computation and inversion. The method can be naturally extended to solve composite optimization problems with non-smooth regularization.
- We theoretically show under certain conditions the proposed approach can improve state-of-the-art iteration complexity, and empirically demonstrate it can substantially improve the convergence of existing methods.

### 1.3 Notations

We use $[n]$ to denote the set $1, ..., n$. For a vector $w \in \mathbb{R}^d$, we use $\|w\|$ to denote its $\ell_2$ norm and $\|w\|_1$ to denote its $\ell_1$ norm. For a matrix $X$, we use $\|X\|_2$ to denote its spectral norm and $\|X\|_F$ to denote its Frobenius norm, and $\lambda_{\min}(X)$ to denote its minimum singular value, and $\mathrm{tr}(X)$ to denote the trace of $X$. We use $I$ to denote an identity matrix. For two symmetric matrices $A$ and $B$, we denote $A \succeq B$ if $A - B$ is positive semi-definite. For two sequences of numbers $\{a_n\}$ and $\{b_n\}$, we say $a_n = \mathcal{O}(b_n)$ if $a_n \leq Cb_n$ for $n$ large enough, with some positive constant $C$. We use the notation $\widetilde{\mathcal{O}}(\cdot)$ to hide poly-log factors. We also use $a_n \lesssim b_n$ to denote $a_n = \mathcal{O}(b_n)$, $a_n \gtrsim b_n$ to denote $b_n = \mathcal{O}(a_n)$, $a_n \asymp b_n$ to denote $a_n = \mathcal{O}(b_n)$ and $b_n = \mathcal{O}(a_n)$.

### 1.4 Organization

The rest of the paper is organized as follows: we introduce the proposed method in Section 2, and present the main theoretical results in Section 3. In Section 4 we present the convergence analysis of inexact minibatch accelerated SVRG update as a key step in proving the main results, and this analysis might be of independent interest. We provide some empirical comparisons over existing approaches in Section 5, and conclusions in Section 6. Some detailed proofs are deferred to the Appendices.

## 2 Minibatch Stochastic Variance Reduced Proximal Iterations

In this section we present the proposed approach for minimizing (1), which naturally integrate both noisy first-order and higher-order information in an efficient way.

### 2.1 Algorithm Description

The high-level description is given in Algorithm 1, which consists of two main innovations:

- simultaneous incorporation of first-order and higher-order information via sub-sampling;
- allowing larger minibatch sizes through acceleration with inexact minimization.

We explain these features in the remaining of this section.

**Algorithm 1** MB-SVRP: Minibatch Stochastic Variance Reduced Proximal Iterations.

**Parameters** $\eta$, $b$, $\widetilde{\lambda}$, $\nu$, $\varepsilon$.
**Initialize** $\widetilde{w}_0 = 0$.
**Sampling** Sampling $b$ items from $[n]$ to form a minibatch $\bar{B}$.
**for** $s = 1, 2, \ldots$ **do**
  **Calculate** $\widetilde{v} = \frac{1}{n}\sum_{i=1}^{n} \nabla f_i(\widetilde{w}_{s-1})$.
  **Initialize** $y_0 = w_0 = \widetilde{w}_s$.
  **for** $t = 1, 2, \ldots, m$ **do**
    **Sampling** $b$ items from $[n]$ to form a minibatch $B_t$.
    **Find** $w_t$ **that approximately solve** (4) (**Option-I**) or (5) (**Option-II**) such that $\widetilde{f}_t(w_t) - \min_w \widetilde{f}_t(w) \le \varepsilon$.

$$\textbf{Option-I}: \widetilde{f}_t(w) := \frac{1}{b}\sum_{i\in\bar{B}} f_i(w) - \left\langle \frac{1}{b}\sum_{i\in\bar{B}} \nabla f_i(y_{t-1}), w \right\rangle + \left\langle \frac{\eta}{b}\sum_{i\in B_t} \nabla f_i(y_{t-1}) + \eta\widetilde{v} - \frac{\eta}{b}\sum_{i\in B_t} \nabla f_i(\widetilde{w}_{s-1}), w \right\rangle + \frac{\widetilde{\lambda}}{2}\left\|w - y_{t-1}\right\|^2, \quad (4)$$

$$\textbf{Option-II}: \widetilde{f}_t(w) := \frac{1}{2}(w - y_{t-1})^\top \left(\frac{1}{b}\sum_{i\in\bar{B}} \nabla^2 f_i(y_{t-1})\right)(w - y_{t-1}) + \left\langle \frac{\eta}{b}\sum_{i\in B_t} \nabla f_i(y_{t-1}) + \eta\widetilde{v} - \frac{\eta}{b}\sum_{i\in B_t} \nabla f_i(\widetilde{w}_{s-1}), w \right\rangle + \frac{\widetilde{\lambda}}{2}\left\|w - y_{t-1}\right\|^2. \quad (5)$$

    **Update**:

$$y_t = w_t + \nu(w_t - w_{t-1}).$$

  **end for**
  **Update** $\widetilde{w}_s = w_m$.
**end for**
**Return** $\widetilde{w}_s$

### Sampling both First-order and Higher-order Information

We first discuss the main building block of the algorithm. Suppose at iteration $t$, given the previous iterate $w_{t-1}$, we consider the following update rule:

$$w_t = \arg\min_w \frac{1}{b}\sum_{i\in\bar{B}} f_i(w) - \left\langle \frac{1}{b}\sum_{i\in\bar{B}} \nabla f_i(w_{t-1}), w \right\rangle + \left\langle \frac{\eta}{b}\sum_{i\in B_t} \nabla f_i(w_{t-1}) + \eta\nabla f(\widetilde{w}) - \frac{\eta}{b}\sum_{i\in B_t} \nabla f_i(\widetilde{w}), w \right\rangle + \frac{\widetilde{\lambda}}{2}\left\|w - w_{t-1}\right\|^2, \quad (6)$$

where $\bar{B}, B_t$ are some randomly sampled minibatch from $1, ..., n$, both with minibatch size $b$; $\eta$ and $\widetilde{\lambda}$ are stepsize parameters, and $\widetilde{w}$ is a "reference" predictor used for reducing variance. The main feature of updating rule (6) is it considered both noisy first-order and higher-order information via minibatch sampling. The term $\frac{1}{b}\sum_{i\in\bar{B}} f_i(w) - \left\langle \frac{1}{b}\sum_{i\in\bar{B}} \nabla f_i(w_{t-1}), w \right\rangle$ in (6) incorporates noisy higher-order information of the objective function. In particular, second-order methods (e.g.

Newton methods) often use the following second-order Taylor approximation of the function:

$$\frac{1}{b}\sum_{i\in\bar{B}} f_i(w) \approx \frac{1}{b}\sum_{i\in\bar{B}} f_i(w_{t-1}) + \left\langle \frac{1}{b}\sum_{i\in\bar{B}} \nabla f_i(w_{t-1}), w - w_{t-1} \right\rangle$$
$$+ \frac{1}{2}(w - w_{t-1})^\top \left( \frac{1}{b}\sum_{i\in\bar{B}} \nabla^2 f_i(w_{t-1}) \right)(w - w_{t-1}) . \tag{7}$$

If we plug the second-order approximation (7) into (6), we get the update rule as

$$w_t \approx \arg\min_w \frac{1}{2}(w - w_{t-1})^\top \left( \frac{1}{b}\sum_{i\in\bar{B}} \nabla^2 f_i(w_{t-1}) \right)(w - w_{t-1}) + \frac{\widetilde{\lambda}}{2}\|w - w_{t-1}\|^2$$
$$+ \left\langle \frac{\eta}{b}\sum_{i\in B_t} \nabla f_i(w_{t-1}) + \eta \nabla f(\widetilde{w}) - \frac{\eta}{b}\sum_{i\in B_t}\nabla f_i(\widetilde{w}), w \right\rangle \tag{8}$$
$$= w_{t-1} - \eta \left( \frac{1}{b}\sum_{i\in\bar{B}} \nabla^2 f_i(w_{t-1}) + \widetilde{\lambda} I \right)^{-1} \left( \frac{1}{b}\sum_{i\in B_t} \nabla f_i(w_{t-1}) + \nabla f(\widetilde{w}) - \frac{1}{b}\sum_{i\in B_t} \nabla f_i(\widetilde{w}) \right), \tag{9}$$

which can be treated as a variant of sub-sampled Newton method, combined with the minibatch stochastic gradient with variance reduction. Moreover, when $f_i(w)$ is a quadratic function of $w$, the approximation (7) is exact (thus Option-I and Option-II in Algorithm 1 are coincident). Therefore the update rule (6) can be treated as a preconditioned minibatch SVRG update rule, with $\left(\frac{1}{b}\sum_{i\in\bar{B}} \nabla^2 f_i(w_{t-1}) + \widetilde{\lambda} I\right)^{-1}$ as a precondition matrix. This is formalized in the proposition below.

**Proposition 2.** *When $f_i(w), \forall i$ is a quadratic function of $w$, then the update rule of (6) is equivalent to the following preconditioned minibatch SVRG update rule:*

$$w_t \leftarrow w_{t-1} - \eta \left( \widetilde{H} + \widetilde{\lambda} I \right)^{-1} \left( \frac{1}{b}\sum_{i\in B_t} \nabla f_i(w_{t-1}) + \nabla f(\widetilde{w}) - \frac{1}{b}\sum_{i\in B_t} \nabla f_i(\widetilde{w}) \right),$$

*where*

$$\widetilde{H} = \frac{1}{b}\sum_{i\in\bar{B}} \nabla^2 f_i(w)$$

*is the sub-sampled Hessian matrix.*

### Connection to Existing Methods

Depending on the choice of these parameters, we observe that the update rule of (6) can be viewed as a generalization of several update rules proposed recently:

- When $b = 1$, $\widetilde{\lambda} \to \infty$ and $\eta = \frac{\widetilde{\lambda}}{L}$, the term $\frac{1}{b}\sum_{i\in\bar{B}} f_i(w) - \langle \frac{1}{b}\sum_{i\in\bar{B}} \nabla f_i(w_{t-1}), w\rangle$ is negligible, thus (6) reduced to standard SVRG update (Johnson and Zhang, 2013):

$$w_t \leftarrow w_{t-1} - \frac{1}{L}\left( \nabla f_i(w_{t-1}) + \frac{1}{n}\sum_{i=1}^n \nabla f_i(\widetilde{w}) - \nabla f_i(\widetilde{w}) \right).$$

- When $b > 1$, $\widetilde{\lambda}, \eta \to \infty$, then the term $\frac{1}{b}\sum_{i\in\bar{B}} f_i(w) - \left\langle \frac{1}{b}\sum_{i\in\bar{B}} \nabla f_i(w_{t-1}), w\right\rangle$ in (6) is negligible as well, (6) reduces to the following update

$$w_t \leftarrow w_{t-1} - \frac{\eta}{\widetilde{\lambda}}\left(\frac{1}{b}\sum_{i\in B_t}\nabla f_i(w_{t-1}) + \frac{1}{n}\sum_{i=1}^{n}\nabla f_i(\widetilde{w}) - \frac{1}{b}\sum_{i\in B_t}\nabla f_i(\widetilde{w})\right),$$

which recovers the update rule of minibatch semi-stochastic gradient methods (a.k.a minibatch SVRG) (Konečný et al., 2016).

- When $\bar{B} = B_t$, $b = 1$, $\eta = 1$, (6) reduced to stochastic variance reduced proximal iterations. More specifically, (6) will be reduced to

$$w_t = \arg\min_w f_i(w) + \left\langle \frac{1}{n}\sum_{i=1}^{n}\nabla f_i(\widetilde{w}) - \nabla f_i(\widetilde{w}), w\right\rangle + \frac{\widetilde{\lambda}}{2}\left\|w - w_{t-1}\right\|^2, \tag{10}$$

by checking the first order optimality condition it is easy to verify that (10) is performing the following update:

$$w_t \leftarrow w_{t-1} - \frac{1}{\widetilde{\lambda}}\left(\nabla f_i(w_t) + \frac{1}{n}\sum_{i=1}^{n}\nabla f_i(\widetilde{w}) - \nabla f_i(\widetilde{w})\right),$$

compared with standard SVRG update with stepsize $\frac{1}{\widetilde{\lambda}}$, we see the only difference is the gradient evaluation on $f_i$ is based on "future" iterate $w_t$ rather than the current iterate $w_{t-1}$ in SVRG. Stochastic proximal iterations based on SAGA method has been analyzed in (Defazio, 2016) recently, and shown to achieve accelerated rate without explicit momentum step.

- When $B_t = \bar{B}$, $b > 1$, $\eta = 1$, and $\widetilde{w} = w_{t-1}$, (6) will be reduced to

$$w_t = \arg\min_w \frac{1}{b}\sum_{i\in\bar{B}} f_i(w) + \left\langle \frac{1}{n}\sum_{i=1}^{n}\nabla f_i(w_{t-1}) - \frac{1}{b}\sum_{i\in\bar{B}}\nabla f_i(w_{t-1}), w\right\rangle + \frac{\widetilde{\lambda}}{2}\left\|w - w_{t-1}\right\|^2, \tag{11}$$

(11) covers the update rule of DANE algorithm (Shamir et al., 2014), which is a communication-efficient distributed optimization algorithm. DANE uses the data on local machine to form the minibatch $\bar{B}$ and every round machines communicate the gradient vector based on their local data. As shown in (Shamir et al., 2014), DANE is provably communication more efficient than the first-order methods in certain scenarios.

- When $\bar{B} = B_t$, and we ignore the linear term $-\left\langle \frac{1}{b}\sum_{i\in\bar{B}} \nabla f_i(w_{t-1}), w\right\rangle$ $+\left\langle \frac{\eta}{b}\sum_{i\in B_t}\nabla f_i(w_{t-1}) + \eta\nabla f(\widetilde{w}) - \frac{\eta}{b}\sum_{i\in B_t}\nabla f_i(\widetilde{w}), w\right\rangle$, then (6) reduces to the minibatch proximal iterations (Li et al., 2014; Wang et al., 2017a):

$$w_t = \arg\min_w \frac{1}{b}\sum_{i\in B_t} f_i(w) + \frac{\widetilde{\lambda}}{2}\left\|w - w_{t-1}\right\|^2.$$

Such an update allows larger minibatch size than standard minibatch SGD, but without the linear correction term as we considered in (6). For such methods only sublinear convergence can be established for finite-sum problems.

**Large Minibatch Size via Acceleration with Inexact Minimization**

Using (6) as the building block, we propose the MB-SVRP (minibatch stochastic variance reduced proximal iterations) method, which is detailed in Algorithm 1. At the beginning of the algorithm, we form a minibatch $\bar{B}$ by sampling from $1, ..., n$ and fix it for the whole optimization process[1]. Then following the SVRG method (Johnson and Zhang, 2013), the algorithm is divided to multiple stages, indexed by $s$. At each stage, we iteratively solve a minimization problem of the form (4) (Option I) or (5) (Option II) based on the randomly sampled minibatch $B_t$.

Compared with the simple update rule in (6), the major difference in Algorithm 1 is that we consider a momentum scheme by maintaining two sequences $\{w_t, y_t\}$, which is inspired by Nesterov's acceleration technique (Nesterov, 2004) and its recent SVRG variant (Nitanda, 2014). The main theoretical advantage over minibatch SVRG without momentum (Konečný et al., 2016) is that such an acceleration allows us to use a much larger minibatch size (up to a size of $O(\sqrt{n})$) without slowing down the convergence. In comparison the iteration complexity of standard minibatch SVRG will become worse when the minibatch size (Konečný et al., 2016) increases. The advantage of allowing larger minibatch size using acceleration in SVRG type algorithms is in analogy to the situation of stochastic gradient descent (without variance reduction) type algorithms (Dekel et al., 2012; Cotter et al., 2011; Lan, 2012).

Since it is often expensive to find the exact minimizer of (4), we consider an approximate minimizer with objective suboptimality $\varepsilon$. When we choose the appropriate $\widetilde{\lambda}$ for $\widetilde{f}_t(w)$ to obtain enough strong convexity, the subproblems defined in Option-I and Option-II can be both be solved to high accuracy efficiently. For example, using SVRG-type algorithms (Johnson and Zhang, 2013) to solve (4) or (5) allows us to find an approximate solution with a small suboptimality $\varepsilon$ with a constant number of passes over the data in the minibatch $\bar{B} \cup B_t$, if $\widetilde{\lambda}$ is set appropriately (discussed in details in our theoretical analysis). This approach avoids Hessian matrix construction and inversion operations, which are often computationally expensive for large-scale problems. Allowing error in gradient oracle has been analyzed in several first-order methods (Schmidt et al., 2011; Devolder et al., 2014) in the batch setting, but has been largely unexplored in stochastic gradient methods. The recent work of (Wang et al., 2017a) analyzed inexact minibatch proximal updates, which can be treated as an implicit minibatch stochastic gradient with errors. In Section 4, we establish convergence rate of *inexact*, *minibatched*, *accelerated* SVRG method, which might be of independent interest.

## 2.2 Extension to Composite Minimization

For many methods that try to incorporate second-order information such as sub-sampled Newton and L-BFGS type algorithms, it is not clear how to extend them to solving non-smooth composite problems. In contrast, the proposed approach can be easily extended to handle non-smooth composite minimization problems. Consider the minimization of:

$$F(w) := \frac{1}{n} \sum_{i=1}^{n} f_i(w) + g(w), \tag{12}$$

where the component functions $f_i(w), \forall i \in [n]$ are smooth and strongly convex, and $g(w)$ is a non-smooth regularizer. For example when $g(w) = \mu \|w\|_1$ and $f_i(w) = \frac{1}{2}(b_i - \langle w, x_i \rangle)^2$, we get the Lasso objective (Tibshirani, 1996); when $g(w) = \mu \|w\|_1 + \frac{\lambda}{2} \|w\|^2$ and $f_i(w) = \log(1 + \exp(-b_i\langle w, x_i\rangle))$, we get the elastic net regularized logistic regression (Zou and Hastie, 2005).

[1] In practice, we can also consider the varying $\bar{B}$ option by simply setting $\bar{B} = B_t$, which makes the algorithm simpler to implement. In practice, we observe no significant difference between this option and the pre-fixed $\bar{B}$. Here we consider fixed $\bar{B}$ mainly for the sake of simplicity in our theoretical analysis.

We can easily modify Algorithm 1 to solve (12). The idea is rather straightforward: at each inner iteration, rather than (4), we simply solve the following minibatch composite minimization problem:

$$w_t \approx \arg\min_w \widetilde{F}_t(w) := \frac{1}{\eta b}\sum_{i\in\bar{B}} f_i(w) - \left\langle \frac{1}{\eta b}\sum_{i\in\bar{B}} \nabla f_i(y_{t-1}), w \right\rangle$$
$$+ \left\langle \frac{1}{b}\sum_{i\in B_t} \nabla f_i(y_{t-1}) + \widetilde{v} - \frac{1}{b}\sum_{i\in B_t} \nabla f_i(\widetilde{w}_{s-1}), w \right\rangle + \frac{\widetilde{\lambda}}{2\eta}\|w - y_{t-1}\|^2 + g(w). \quad (13)$$

Here we slightly re-scale the approximated loss term to make sure that the relative weight between the approximated loss and the regularization term is correct. Since the objective (13) is a standard finite-sum composite minimization problem we could apply prox-SVRG or prox-SAGA to solve (13) efficiently when the term $\frac{\widetilde{\lambda}}{2\eta}\|w - y_{t-1}\|^2$ gives sufficient strong convexity.

## 3 Convergence Analysis

In this section we present theoretical results for the proposed approach. For general convex objectives, our convergence analysis is local, i.e. we assume the initial and subsequent solutions are within a bounded region near the optimum, which is common in the analysis of second-order methods in optimization (Nesterov, 2004; Boyd and Vandenberghe, 2004). As will be discussed in the sequel, the size of the region depends on the condition number of the problem, as well as the Lipschitz parameter of the Hessian. For quadratic objectives, the radius of the region can be infinite, which implies global fast convergence.

In this section, for simplicity, we focus on the analysis of Option-II in Algorithm 1 which uses second-order approximation to construct each sub-problem (5). Similar results hold for Option-I because locally around the optimal solution, the two methods are equivalent.

Besides the strong convexity and smoothness, we also need the following Hessian Lipschitz condition.

**Condition 1.** *For each component function $f_i(w), \forall i \in [n]$, its Hessian matrix is $M$-Lipschitz, i.e.*

$$\left\|\nabla^2 f_i(w) - \nabla^2 f_i(w')\right\| \leq M \left\|w - w'\right\|, \forall w, w' \in \mathbb{R}^d, \forall i \in [n].$$

We consider the component function in form of convex loss function of linear predictor: $f_i(w) = \ell(w^\top x_i; b_i)$, for which we have the following formula of Hessian matrix:

$$H_\lambda(w) = \frac{1}{n}\sum_{i=1}^n \nabla^2 f_i(w) = \frac{1}{n}\sum_{i=1}^n \ell''(w^\top x_i; b_i) x_i x_i^\top + \lambda I$$

be the Hessian matrix at $w$, and suppose we have an approximated Hessian

$$\widetilde{H}_{\widetilde{\lambda}}(w_{t-1}) = \frac{1}{b}\sum_{i\in\bar{B}} \nabla^2 f_i(w_{t-1}) + \widetilde{\lambda} I,$$

the following proposition states that performing the preconditioned minibatch SVRG update is equivalent to the standard minibatch SVRG in a linear transformed (preconditioned) space:

**Proposition 3.** *Considering the following minibatch SVRG type update when solving* $\min_w f(w)$*:*

$$w_t \leftarrow w_t - \eta \widetilde{H}_{\widetilde{\lambda}}(w_{t-1})^{-1} \left( \frac{1}{b} \sum_{i \in B_t} \nabla f_i(w_{t-1}) + \frac{1}{n} \sum_{i=1}^{n} \nabla f_i(\widetilde{w}) - \frac{1}{b} \sum_{i \in B_t} \nabla f_i(\widetilde{w}) \right), \tag{14}$$

*which is equivalent to solving a minimization problem with respect to* $z$*:*

$$\min_z f(\widetilde{H}_{\widetilde{\lambda}}^{-1/2}(w^*)z)$$

*via the variable transform* $z = \widetilde{H}_{\widetilde{\lambda}}^{1/2}(w^*)w$*. Moreover, the update rule* (14) *is equivalent to the following update on* $z$*:*

$$\begin{aligned} z_t \leftarrow z_{t-1} - \eta &\left( \frac{1}{b} \sum_{i \in B_t} \nabla_z f_i(\widetilde{H}_{\widetilde{\lambda}}^{-1/2}(w^*)z_{t-1}) + \frac{1}{n} \sum_{i=1}^{n} \nabla_z f_i(\widetilde{H}_{\widetilde{\lambda}}^{-1/2}(w^*)\widetilde{z}) - \frac{1}{b} \sum_{i \in B_t} \nabla_z f_i(\widetilde{H}_{\widetilde{\lambda}}^{-1/2}(w^*)\widetilde{z}) \right) \\ &+ \eta \left( \widetilde{H}_{\widetilde{\lambda}}(w^*)^{-1} - \widetilde{H}_{\widetilde{\lambda}}(w_{t-1})^{-1} \right) \left( \frac{1}{b} \sum_{i \in B_t} \nabla f_i(w_{t-1}) + \frac{1}{n} \sum_{i=1}^{n} \nabla f_i(\widetilde{w}) - \frac{1}{b} \sum_{i \in B_t} \nabla f_i(\widetilde{w}) \right). \end{aligned} \tag{15}$$

With the above proposition, we can see that the proposed update (8) is implicitly performing inexact minibatch SVRG update on the following transformed problem:

$$\min_z f(\widetilde{H}_{\widetilde{\lambda}}^{-1/2}(w^*)z) := \frac{1}{n} \sum_{i=1}^{n} \ell(w^\top \widetilde{H}_{\widetilde{\lambda}}^{-1/2}(w^*)z; b_i) + \lambda \left\| \widetilde{H}_{\widetilde{\lambda}}^{-1/2}(w^*)w \right\|^2. \tag{16}$$

As a straightforward extension, it is also easy to see that using the acceleration technique as Algorithm 1 is equivalent to performing inexact *accelerated* minibatch SVRG update, similar to the exact update analyzed in (Nitanda, 2014).

The benefit of the transformed problem is that when $\widetilde{H}_{\widetilde{\lambda}}(w^*) \approx H_\lambda(w^*)$, the condition number of the new problem (16) becomes smaller than that of the original problem. As a consequence, when the computational cost of update (14) is cheap, we expect improved convergence and runtime guarantees for solving the original problem. Before analyzing the condition number, we state the following lemma, which connects the inexactness in $w$ space when performing the update (8) (and (14)) to the inexactness in $z$ space when performing the update (15).

**Lemma 4.** *Let* $\widetilde{f}_t(w)$ *be the function to be minimized in* (8)*, and* $w_t^* = \arg\min_w \widetilde{f}_t(w)$ *to be its exact minimizer. Let*

$$\bar{z}_t = z_{t-1} - \eta \left( \frac{1}{b} \sum_{i \in B_t} \nabla_z f_i(\widetilde{H}_{\widetilde{\lambda}}^{-1/2}(w^*)z_{t-1}) + \frac{1}{n} \sum_{i=1}^{n} \nabla_z f_i(\widetilde{H}_{\widetilde{\lambda}}^{-1/2}(w^*)\widetilde{z}) - \frac{1}{b} \sum_{i \in B_t} \nabla_z f_i(\widetilde{H}_{\widetilde{\lambda}}^{-1/2}(w^*)\widetilde{z}) \right),$$

*be the exact minibatch SVRG update solution in the preconditioned space, and* $z_t = \widetilde{H}_{\widetilde{\lambda}}^{1/2}(w^*)w_t$ *be the actual update in precondition space by approximately solving* (8)*. Let* $z_t^* = \widetilde{H}_{\widetilde{\lambda}}^{1/2}(w^*)w_t^*$ *be the exact minimizer in the preconditioned space, and let* $\vartheta_t = \frac{1}{b}\sum_{i \in B_t} \nabla f_i(w_{t-1}) + \frac{1}{n}\sum_{i=1}^n \nabla f_i(\widetilde{w}) - \frac{1}{b}\sum_{i \in B_t} \nabla f_i(\widetilde{w})$*. We have the following upper bound on* $\|z_t - \bar{z}_t\|$*:*

$$\|z_t - \bar{z}_t\| \le \sqrt{\frac{2(\widetilde{f}_t(w_t) - \widetilde{f}_t(w_t^*))(L + \widetilde{\lambda})}{\lambda + \widetilde{\lambda}}} + \frac{\eta M \|w_{t-1} - w^*\| \|\vartheta_t\|}{\widetilde{\lambda}^2}.$$

**Remark 1.** *Lemma 4 states that the inexactness in the minibatch SVRG update in the $z$ space can be decomposed to two sources, the first source is the fact that we solve the minimization problem of objective $\widetilde{f}_t(w)$ inexactly, and the second source is the fact that the Hessian matrix is changing at different locations. Nevertheless, we can fix the $z$ space that use Hessian matrix at optimum $\widetilde{H}_{\widetilde{\lambda}}(w^*)$ as the preconditioner. For quadratic objectives where $M = 0$, the second error term becomes zero.*

## 3.1 Improved Condition Number

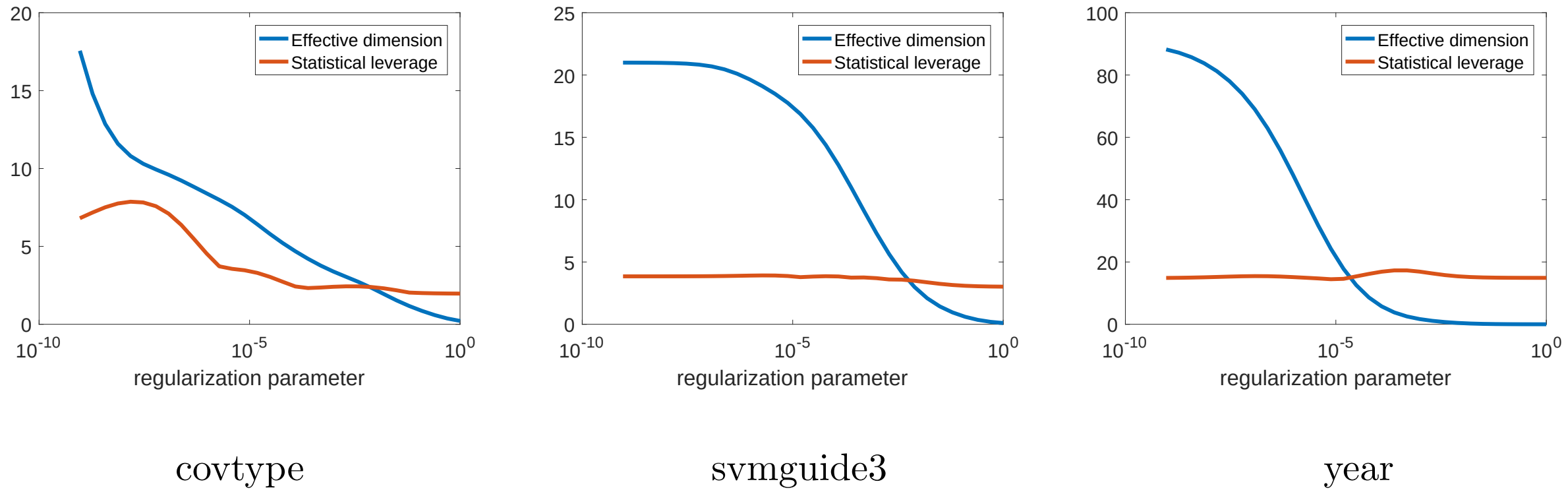

Figure 1: Effective dimension and maximum statistical leverage with different regularization parameters $\widetilde{\lambda}$ on some empirical datasets for ridge regression problems.

In this section we analyze the condition number in the "preconditioned" space (16), where $\bar{B}$ is sampled uniformly from $1, ..., n$. For the analysis, we introduce two notions which describes the global and local properties of data. The following definitions of effective dimension and bounded statistical leverage extend the ones defined in (Hsu et al., 2014) which were used to analyzing the generalization performance of ridge regression.

**Definition 5.** *(**Effective dimension**) Let the $\lambda_1, ..., \lambda_d$ be the top-$d$ eigenvalues of $H_0(w^*) = (1/n)\sum_{i=1}^n \ell''(w^{*\top} x_i; b_i) x_i x_i^\top$, define the effective dimension $d_{\widetilde{\lambda}}$ (for some $\widetilde{\lambda} \geq 0$) of $H_0(w^*)$ being*

$$d_{\widetilde{\lambda}} = \sum_{j=1}^{d} \frac{\lambda_j}{\lambda_j + \widetilde{\lambda}}.$$

We see $d_{\widetilde{\lambda}}$ is a decreasing function of $\widetilde{\lambda}$, and $d_0 = d$ when $\widetilde{\lambda} = 0$. When the spectrum of $\frac{1}{n}\sum_{i=1}^n \ell''(w^{*\top} x_i; b_i) x_i x_i^\top$ is decaying very fast, the effective dimension $d_{\widetilde{\lambda}}$ can be significantly smaller than $d$ for moderate $\widetilde{\lambda}$. The following notion of statistical leverage have been used in regression analysis (Chatterjee and Hadi, 2009; Hsu et al., 2014) and matrix approximation (Mahoney et al., 2011).

**Definition 6.** *(**Statistical leverage at $\widetilde{\lambda}$**) Let $H_{\widetilde{\lambda}}(w^*) = \frac{1}{n}\sum_{i=1}^n \ell''(w^{*\top} x_i; b_i) x_i x_i^\top + \widetilde{\lambda} I$, we say the statistical leverage of data matrix $X$ is bounded by $\rho_{\widetilde{\lambda}}$ at $\widetilde{\lambda}$ if*

$$\max_{i \in [n]} \frac{\left\| H_{\widetilde{\lambda}}^{-1/2}(w^*) \ell''(w^{*\top} x_i; b_i)^{1/2} x_i \right\|}{\sqrt{(1/n)\sum_{j=1}^n \left\| H_{\widetilde{\lambda}}^{-1/2}(w^*) \ell''(w^{*\top} x_j; b_j)^{1/2} x_j \right\|^2}} \leq \rho_{\widetilde{\lambda}}. \tag{17}$$

For ridge regression problems, the above definition is slightly different from the one used in (Hsu et al., 2014) in the sense that the empirical Hessian matrix $H_{\widetilde{\lambda}}$ in (17) is replaced by the population Hessian matrix $\mathbb{E}\left[x_i^\top x_i\right] + \widetilde{\lambda} I$ in (Hsu et al., 2014). When the sample size $n$ is large, the differences in these two definitions are minor. As argued in (Hsu et al., 2014), when $x_i$ is drawn from sub-Gaussian distributions, $\left(\mathbb{E}\left[x_i^\top x_i\right]\right)^{-1/2} x_i$ is isotropic. In this case the statistical leverage only grows logarithmically with dimension. Our theoretical analysis relies on the condition that the two quantities defined above are not too large. This condition often holds in practice, as shown in Figure 1, which plots the effective dimension and the statistical leverage for several real world datasets with varying regularization parameters.

The convergence of stochastic gradient based algorithms for minimizing the objective (40) depends on two important quantities: strong convexity and smoothness. The strong convexity for problem (40) is $\lambda_{\min}\left(\widetilde{H}_{\widetilde{\lambda}}^{-1/2}(w^*)\left(\frac{1}{n}\sum_{i=1}^n \ell''(w^\top x_i; b_i)x_i x_i^\top + \lambda I\right)\widetilde{H}_{\widetilde{\lambda}}^{-1/2}(w^*)\right)$, which is a global property of the objective function. However, we must consider the smoothness parameter for individual function $f_i(\widetilde{H}_{\widetilde{\lambda}}^{-1/2}(w^*)z)$, which in our context, is $\max_i \left\{\ell''(w^{*\top} x_i; b_i)^{1/2} x_i^\top \widetilde{H}_{\widetilde{\lambda}}^{-1}(w^*)x_i\right\}$. Both of these quantities are related to how close is the constructed Hessian approximation $\widetilde{H}_{\widetilde{\lambda}}(w^*)$ to the true Hessian $H_{\widetilde{\lambda}}(w^*)$ in spectral norm. The following lemma bound this quantity using matrix concentration.

**Lemma 7.** *If $\bar{B}$ is formed by uniform sampling with replacement from $[n]$, then we have the following concentration bound, with probability at least $1-\delta$,*

$$\left\| H_{\widetilde{\lambda}}^{-1}(w^*)\left(\frac{1}{b}\sum_{i\in\bar{B}} \ell''(w^{*\top} x_i; b_i)x_i x_i^\top - \frac{1}{n}\sum_{i=1}^n \ell''(w^{*\top} x_i; b_i)x_i x_i^\top\right)\right\|_2 \le 2\log\left(\frac{d}{\delta}\right)\cdot\sqrt{\frac{\rho_{\widetilde{\lambda}}^2 d_{\widetilde{\lambda}}}{b}}.$$

### Strong Convexity

Based on above lemma, we have the following lower bound of the strong convexity for (16), specified in the following lemma.

**Lemma 8.** *If we set $\widetilde{\lambda}$ such that $\widetilde{\lambda} \ge \lambda$, then with probability at least $1-\delta$ over the random choice of $\bar{B}$ to form $\widetilde{H}_{\widetilde{\lambda}}$, we have*

$$\min_{\|w-w^*\|\le R} \lambda_{\min}\left(\widetilde{H}_{\widetilde{\lambda}}^{-1/2}(w^*)H_\lambda(w)\widetilde{H}_{\widetilde{\lambda}}^{-1/2}(w^*)\right) \ge \frac{\lambda}{\widetilde{\lambda}}\frac{1}{1+2\log(d/\delta)\sqrt{(\rho_{\widetilde{\lambda}}^2 d_{\widetilde{\lambda}})/b} + MLR + MR/\widetilde{\lambda}}.$$

### Smoothness

Next we explore the smoothness parameter of (16), which is upper bounded by

$$\max_{i\in[n]} \ell''(w^{*\top} x_i; b_i)x_i^\top \widetilde{H}_{\widetilde{\lambda}}^{-1}(w^*)x_i.$$

The most straightforward way to upper bound $\max_{i\in[n]} \ell''(w^{*\top} x_i; b_i)x_i^\top \widetilde{H}_{\widetilde{\lambda}}^{-1}(w^*)x_i$ is

$$\begin{aligned}\max_{i\in[n]} \ell''(w^{*\top} x_i; b_i)x_i^\top \widetilde{H}_{\widetilde{\lambda}}^{-1}(w^*)x_i &\le \max_{i\in[n]} \ell''(w^{*\top} x_i; b_i)\|x_i\|^2 \lambda_{\min}(\widetilde{H}_{\widetilde{\lambda}}^{-1}(w^*)) \\ &= \max_{i\in[n]} \frac{\ell''(w^{*\top} x_i; b_i)\|x_i\|^2}{\widetilde{\lambda}} = \frac{L}{\widetilde{\lambda}}.\end{aligned}$$

In this way, the condition number after preconditioning becomes

$$\mathcal{O}\left(\frac{4\widetilde{\lambda}}{\lambda} \cdot \frac{L}{\widetilde{\lambda}}\right) = \frac{4L}{\lambda},$$

which didn't show any advantage of using preconditioning. This is not surprising because it is known that for high dimensional problems, the worse case behavior of second order methods is no better than that of the first order methods.

In the lemma below, we provide an improved analysis of smoothness which is based on the notion of effective dimension (Definition 5) and statistical leverage (Definition 6).

**Lemma 9.** *If we choose* $\widetilde{\lambda} \geq \lambda$ *and* $b \geq 16\rho_{\widetilde{\lambda}}^2 d_{\widetilde{\lambda}} \log^2(d/\delta)$*, then with probability at least* $1-\delta$ *over the random choice of* $\bar{B}$ *to form* $\widetilde{H}_{\widetilde{\lambda}}(w^*)$*, we have*

$$\max_{i\in[n]} \ell''({w^*}^\top x_i; b_i) x_i^\top \widetilde{H}_{\widetilde{\lambda}}^{-1}(w^*) x_i \leq 2\rho_{\widetilde{\lambda}}^2 d_{\widetilde{\lambda}}.$$

Combining Lemma 8 and Lemma 9 we get the following result about the condition number for (16). If we choose $\widetilde{\lambda}$ and $b$ that satisfy

$$\widetilde{\lambda} \geq \max\left\{\lambda, \frac{L}{b}\right\}, \quad b \geq 16\rho_{\widetilde{\lambda}}^2 d_{\widetilde{\lambda}} \log^2\left(\frac{d}{\delta}\right), \quad R \leq \min\left\{\frac{1}{4ML}, \frac{\widetilde{\lambda}}{4M}\right\},$$

then with probability at least $1-\delta$, the condition number for stochastic gradient algorithms after preconditioning scales as $\frac{4\rho_{\widetilde{\lambda}}^2 d_{\widetilde{\lambda}} \widetilde{\lambda}}{\lambda}$, more specifically when $\widetilde{\lambda} = \max\left\{\lambda, \frac{L}{b}\right\}$, then the condition number of (16) can be upper bounded by:

$$\max\left\{4\rho_{\widetilde{\lambda}}^2 d_{\widetilde{\lambda}}, \frac{4L\rho_{\widetilde{\lambda}}^2 d_{\widetilde{\lambda}}}{\lambda b}\right\},$$

which improves the original condition number $L/\lambda$ by a factor of at least $\widetilde{\mathcal{O}}(\rho_{\widetilde{\lambda}}^2 d_{\widetilde{\lambda}}/b)$.

## 3.2 Improved Runtime Guarantee

Based on the above analysis, we have the following main results stating the iteration complexity of MP-SVRG algorithms on general convex objectives (1).

**Theorem 10.** *Consider Algorithm 1 on problems with* $L$*-smooth and* $\lambda$*-strongly convex functions that satisfied Condition 1 with Hessian Lipschitz parameter* $M$*. Suppose we sample minibatch* $\bar{B}$ *uniformly from* $1, ..., n$*, and set the tuning parameters as*

$$\widetilde{\lambda} = \max\left\{\lambda, \frac{L}{b}\right\}, \quad b \asymp \min\left\{n, \rho_{\widetilde{\lambda}}^2 d_{\widetilde{\lambda}} \left(\frac{L}{\lambda}\right)^{1/3}\right\}, \quad \eta \asymp \min\left\{\frac{b^3\lambda}{(\rho_{\widetilde{\lambda}}^2 d_{\widetilde{\lambda}})^2 L}, \frac{1}{\rho_{\widetilde{\lambda}}^2 d_{\widetilde{\lambda}}}\right\},$$

*and suppose we start from the initialization point that is sufficiently close to the optimum:*

$$\|w_0 - w^*\| \leq R = \mathcal{O}\left(\frac{\lambda^4}{L^2 M}\right). \tag{18}$$

*Given $\epsilon > 0$, and let*

$$\varepsilon \leq \frac{1}{10^5} \cdot \left(\frac{\lambda}{L}\right)^7 \epsilon. \tag{19}$$

*Then for the MP-SVRP algorithm to find the approximate solution $\widetilde{w}_s$ of (1) that reaches the expected $\epsilon$-objective suboptimality*

$$\mathbb{E}f(\widetilde{w}_s) - \min_w f(w) \leq \epsilon,$$

*the total number of oracle calls to solve (8) is no more than*

$$\mathcal{O}\left(\max\left\{\left(\frac{L}{\lambda}\right)^{1/3}, \frac{\rho_{\widetilde{\lambda}}^2 d_{\widetilde{\lambda}} L}{\lambda n^2}\right\} \cdot \log\left(\frac{1}{\epsilon}\right)\right).$$

Note that in Algorithm 1, we can use SVRG to solve each subproblem in (8). This leads to the following result.

**Corollary 11.** *Assume that the conditions of Theorem 10 hold. Moreover, if we use SVRG to solve each subproblem in (8) up to the suboptimality $\varepsilon$, then the total number of gradient evaluations used in the whole MP-SVRP algorithm can be upper bounded by*

$$\mathcal{O}\left(\max\left\{\rho_{\widetilde{\lambda}}^2 d_{\widetilde{\lambda}} \left(\frac{L}{\lambda}\right)^{2/3}, \frac{\rho_{\widetilde{\lambda}}^2 d_{\widetilde{\lambda}} L}{\lambda n}\right\} \cdot \log\left(\frac{L}{\lambda}\right) \cdot \log^2\left(\frac{1}{\epsilon}\right) + n \cdot \log\left(\frac{1}{\epsilon}\right)\right).$$

We make the following remarks about Theorem 10 and Corollary 11.

**Remark 2.** *The solution accuracy requirement in (19) is conservative, but does not affect the overall iteration complexity because solving (8) to high precision only introduces an additional logarithmic factor using SVRG. In practice, we observe that using only one pass SVRG warm started from previous iterate performs well (please refer to the experiments).*

**Remark 3.** *For quadratic objectives, where $M = 0$, the local fast region defined in (18) is the whole space, implying global fast convergence for quadratic objectives. Moreover, in experiments, we observe that even for the general non-quadratic objectives, the proposed MB-SVRP method also converges fast globally. Therefore the size of the local fast convergence region (18) established here might be too conservative. We leave the investigation on improving the size of the local fast convergence region to future work.*

**Remark 4.** *Corollary 11 states that the iteration complexity of Algorithm 1 is*

$$\widetilde{\mathcal{O}}\left(n + \max\left\{\rho_{\widetilde{\lambda}}^2 d_{\widetilde{\lambda}} \left(\frac{L}{\lambda}\right)^{2/3}, \frac{\rho_{\widetilde{\lambda}}^2 d_{\widetilde{\lambda}} L}{\lambda n}\right\}\right),$$

*where $\widetilde{\mathcal{O}}(\cdot)$ hides the minor poly logarithmic factors. When $\rho_{\widetilde{\lambda}}^2 d_{\widetilde{\lambda}}$ is small so that it can be treated as a constant, then the iteration complexity of MB-SVRP improves over standard SVRG by a factor of $\min\left\{n, \left(\frac{L}{\lambda}\right)^{1/3}\right\}$ when the condition number $\frac{L}{\lambda}$ is much larger than $n$. Even compared with accelerated methods (such as SVRG equipped with catalyst acceleration (Shalev-Shwartz and Zhang, 2016; Frostig et al., 2015; Lin et al., 2015a)), it can sometimes better (when $n < L/\lambda < n^2$, see Table 1) for details.*

**Remark 5.** *MB-SVRP method is able to improve the iteration complexity over SVRG when $\rho_{\widetilde{\lambda}}^2 d_{\widetilde{\lambda}}$ is small, thus it is also possible to use MB-SVRP itself to solve the minibatch subproblem (4). Such a nested approach allows us to choose a even smaller $\widetilde{\lambda}$, thus can further reduce the dependency on condition number in the iteration complexity. However, this is at the cost of increasing the dependency on $\rho_{\widetilde{\lambda}}^2 d_{\widetilde{\lambda}}$ and more complicated implementation. We will leave such investigation to future research.*

### 3.3 With Catalyst Acceleration

**Algorithm 2** Acc-MB-SVRP: Accelerated MB-SVRP method.

**Initialize** $w_0 = z_0 = 0$.
**for** $r = 1, 2, \ldots$ **do**
  **Call MP-SVRP** algorithm 1 to approximately solve

$$w_r \approx \arg\min_w f(w) + \frac{\gamma}{2} \|w - z_{r-1}\|^2 . \tag{20}$$

  **Update**

$$z_r = w_r + \nu_r (w_r - w_{r-1}).$$

**end for**
**Return** $w_r$.

Table 1: Comparison of iteration complexity of various finite-sum quadratic optimization algorithms when effective dimension and statistical leverage are bounded, where we compare the different relative scale of condition number $\kappa = L/\lambda$ and sample size $n$, ignoring logarithmic factors.

| | $\kappa \leq n$ | $\kappa = n^{4/3}$ | $\kappa = n^{3/2}$ | $\kappa = n^2$ | $\kappa = n^3$ |
|---|---|---|---|---|---|
| SVRG | $n$ | $n^{4/3}$ | $n^{3/2}$ | $n^2$ | $n^3$ |
| MB-SVRP | $n$ | $\rho_{\widetilde{\lambda}}^2 d_{\widetilde{\lambda}} \cdot n$ | $\rho_{\widetilde{\lambda}}^2 d_{\widetilde{\lambda}} \cdot n$ | $\rho_{\widetilde{\lambda}}^2 d_{\widetilde{\lambda}} \cdot n^{3/2}$ | $\rho_{\widetilde{\lambda}}^2 d_{\widetilde{\lambda}} \cdot n^{5/2}$ |
| Acc-SVRG | $n$ | $n^{7/6}$ | $n^{5/4}$ | $n^{3/2}$ | $n^2$ |
| Acc-MB-SVRP | $n$ | $\left(\rho_{\widetilde{\lambda}}^2 d_{\widetilde{\lambda}}\right)^{1/2} \cdot n$ | $\left(\rho_{\widetilde{\lambda}}^2 d_{\widetilde{\lambda}}\right)^{1/2} \cdot n$ | $\left(\rho_{\widetilde{\lambda}}^2 d_{\widetilde{\lambda}}\right)^{1/2} \cdot n^{5/4}$ | $\left(\rho_{\widetilde{\lambda}}^2 d_{\widetilde{\lambda}}\right)^{1/2} \cdot n^{7/4}$ |

Corollary 11 stated that if the condition number is not too large: $\frac{L}{\lambda} \leq n^{5/4}$, then MB-SVRP only requires logarithmic passes over data to find a solution with high accuracy. When the condition number is large, MP-SVRP can still be slow. By using the accelerated proximal point framework proposed in (Shalev-Shwartz and Zhang, 2016; Frostig et al., 2015; Lin et al., 2015a), it is possible to obtain an accelerated convergence rate which has milder dependence on condition number, the algorithm is outlined in Algorithm 2, which iteratively, approximately call the original algorithm MB-SVRP to solve an augmented proximal point problem. The main iteration complexity is stated in the theorem below.

**Theorem 12.** *For ill-condition problems where $L/\lambda \geq n^{3/2}$, if we set the parameter $\gamma \asymp \frac{\rho_{\widetilde{\lambda}}^2 d_{\widetilde{\lambda}} L}{n^{3/2}} - \lambda$,*

*Algorithm 2 has iteration complexity of*

$$\widetilde{\mathcal{O}}\left(\left(\rho_{\widetilde{\lambda}}^2 d_{\widetilde{\lambda}}\right)^{1/2} \cdot n^{1/4}\left(\frac{L}{\lambda}\right)^{1/2}\right).$$

**Remark 6.** *From theorem 12 we see the iteration complexity of accelerated MP-SVRP grows at a square root rate with respect to the condition number, which is much better than that of the non-accelerated algorithms, and this is similar to the behavior of accelerated SVRG algorithms. Moreover, accelerated MP-SVRP improves the iteration complexity of accelerated SVRG by a factor of* $n^{1/4}$ *when* $\sqrt{\rho_{\widetilde{\lambda}}^2 d_{\widetilde{\lambda}}}$ *is small.*

# 4 Convergence Analysis of Inexact Accelerated Minibatch SVRG

**Algorithm 3** IMBA-SVRG: Inexact Minibatch Accelerated SVRG Method.

**Parameters** $\alpha = \sqrt{\eta\lambda/2}$.
**Initialize** $\widetilde{w}_0 = 0$.
**for** $s = 1, 2, \ldots$ **do**
  **Calculate** $\widetilde{v} = \frac{1}{n}\sum_{i=1}^n \nabla f_i(\widetilde{w}_{s-1})$.
  **Initialize** $y_0 = w_0 = \widetilde{w}_s$.
  **for** $t = 1, 2, \ldots, m$ **do**
    **Sampling** $b$ items from $[n]$ to form a minibatch $B_t$.
    **Inexact update** $w_t$ such that

$$w_t = \bar{w}_t + \xi_t, \tag{21}$$

    where:

$$\bar{w}_t = y_{t-1} - \eta\left(\frac{1}{b}\sum_{i\in B_t}\nabla f_i(y_{t-1}) - \frac{1}{b}\sum_{i\in B_t}\nabla f_i(\widetilde{w}_{s-1}) + \widetilde{v}\right).$$

    **Update**

$$y_t = w_t + \left(\frac{1-\alpha}{1+\alpha}\right)(w_t - w_{t-1}).$$

  **end for**
  **Update** $\widetilde{w}_s = w_m$.
**end for**
**Return** $\widetilde{w}_s$

The main results established in Section 3 relies on the analysis for inexact, minibatch, accelerated SVRG update (IMBA-SVRG Algorithm 3), where in each stochastic step, we allow a small error $\xi_t$ in the updating. In this section, we show that as long as the inexactness at each iteration is small enough, IMBA-SVRG can use a large minibatch size without slowing down the convergence.

**Theorem 13.** *Consider the IMBA-SVRG algorithm. If we choose the stepsize as*

$$\eta = \min\left\{\frac{b^2\lambda}{6400L^2}, \frac{1}{8L}\right\},$$

*and when the deviation* $\xi_k$ *satisfies* $\forall 1 \le k \le t$

$$\|\xi_k\|^2 \le \frac{2\lambda^2\eta^3\sqrt{\lambda\eta}}{15}(f(\widetilde{w}_{s-1}) - f(w^*)),$$

*and the number of iterations* $t$ *satisfies*

$$t \ge \frac{10}{9\sqrt{\lambda\eta}}\log(36),$$

*then*

$$\mathbb{E}\left[f(\widetilde{w}_s) - f(w^*)\right] \le \frac{1}{2}\left[f(\widetilde{w}_{s-1}) - f(w^*)\right].$$

**Remark 7.** *Depending on the relative magnitude of minibatch size* $b$ *and condition number* $\frac{L}{\lambda}$, *the overall iteration complexity can be divided into the following two situations.*

- $b \le 20\sqrt{\frac{2L}{\lambda}}$. *If we choose the stepsize as* $\eta = \frac{b^2\lambda}{6400L^2}$, *and the deviation* $\xi_k$ *satisfies* $\forall 1 \le k \le t$: $\|\xi_k\|^2 \le \frac{2\lambda^2\eta^3\sqrt{\lambda\eta}}{15}(f(\widetilde{w}_{s-1}) - f(w^*))$, *then* $\mathbb{E}\left[f(\widetilde{w}_s) - f(w^*)\right] \le \frac{1}{2}\left[f(\widetilde{w}_{s-1}) - f(w^*)\right]$ *when* $t \ge \frac{90L}{\lambda b}$.
- $20\sqrt{\frac{2L}{\lambda}} \le b \le n$. *If we choose stepsize as* $\eta = \frac{1}{8L}$, *and when the deviation* $\xi_k$ *satisfies* $\forall 1 \le k \le t$: $\|\xi_k\|^2 \le \frac{2\lambda^2\eta^3\sqrt{\lambda\eta}}{15}(f(\widetilde{w}_{s-1}) - f(w^*))$, *then* $\mathbb{E}\left[f(\widetilde{w}_s) - f(w^*)\right] \le \frac{1}{2}\left[f(\widetilde{w}_{s-1}) - f(w^*)\right]$ *when* $t \ge \sqrt{\frac{8L}{\lambda}}$.

A direct consequence of Theorem 13 is the following iteration complexity of IMBA-SVRG.

**Corollary 14.** *For IMBA-SVRG algorithm, If we choose stepsize as* $\eta = \min\left\{\frac{b^2\lambda}{6400L^2}, \frac{1}{8L}\right\}$, *and when at every state* $s$, *the deviation* $\xi_k$ *at every iteration satisfies* $\forall 1 \le k \le t$: $\|\xi_k\|^2 \le \frac{2\lambda^2\eta^3\sqrt{\lambda\eta}}{15}(f(\widetilde{w}_{s-1}) - f(w^*))$, *then* $\mathcal{O}\left(\log\left(\frac{1}{\epsilon}\right)\right)$ *full gradient evaluations and*

$$\mathcal{O}\left(\max\left\{\frac{L}{\lambda b}, \sqrt{\frac{L}{\lambda}}\right\} \cdot \log\left(\frac{1}{\epsilon}\right)\right)$$

*inexact update steps of form* (21) *are sufficient to ensure* $\mathbb{E}f(\widetilde{w}_s) - f(w^*) \le \epsilon$.

Our convergence proof of IMBA-SVRG relies on the machinery of stochastic estimation sequence (Lin et al., 2014), which originated from the framework of estimation sequence developed in (Nesterov, 2004). The details are given in Appendix B.

The main difference between our analysis and that of (Lin et al., 2014) is that we no longer require the inequalities in estimation sequences hold almost surely. Instead we only make sure they hold in expectation, and thus the quadratic lower bound used to construct the stochastic estimation sequence is also different. IMBA-SVRG can be viewed as an inexact extension of accelerated minibatch SVRG (Nitanda, 2014). Although (Nitanda, 2014) also considered stochastic estimation sequences, here we allow error in the stochastic gradient update, and thus need to construct a different estimation sequence which takes the inexactness into consideration. On the other hand, batch (accelerated) gradient methods with inexact first order oracle have been studied in (Schmidt et al., 2011; Villa et al., 2013; Devolder et al., 2014), and the analysis of IMBA-SVRG extends these results in the context of stochastic gradient methods with variance reduction.

Table 2: List of datasets used in the experiments.

| Name | #Instances | #Features | Task |
|---|---|---|---|
| codrna | 59,535 | 8 | Classification |
| covtype | 581,012 | 54 | Classification |
| svmguide3 | 1,243 | 21 | Classification |
| synthetic-c | 10,000 | 1,000 | Classification |
| cadata | 20,460 | 8 | Regression |
| spacega | 3,107 | 6 | Regression |
| synthetic-r | 20,000 | 2,000 | Regression |
| year | 463,715 | 91 | Regression |

# 5 Experiments

In this section we compare the proposed MP-SVRP algorithm to several state-of-the-art methods for minimizing (1). The datasets are summarized in Table 2, most of which can be download from the LibSVM website[2]. We also considered a synthetic dataset for each task, where the data $\{x_i, b_i\}_{i=1}^n$ are i.i.d. drawn from the following model:

$$\text{Regression}: \quad b_i = \langle x_i, \bar{w} \rangle + a_i, \quad x_i \sim \mathcal{N}(0, \Sigma), \quad a_i \sim \mathcal{N}(0,1), \quad \forall i \in [n],$$

$$\text{Classification}: \quad P(b_i = \pm 1) = \frac{\exp(b_i \langle x_i, \bar{w} \rangle)}{1 + \exp(b_i \langle x_i, \bar{w} \rangle)}, \quad x_i \sim \mathcal{N}(0, \Sigma), \quad \forall i \in [n],$$

where entries of $\bar{w}$ are drawn i.i.d. from $\mathcal{N}(0,1)$. To make the problem ill-conditioned with fast decaying spectrum, we set $\Sigma_{ij} = 2^{-|i-j|/500}, \forall i, j \in [n]$.

## 5.1 Empirical Results for Smooth Optimization

We consider ridge regression and logistic regression models for regression and classification problems, respectively. We normalize the dataset by $x_i \leftarrow x_i / \left( \max_i \|x_i\|^2 \right)$ to ensure that the maximum norm of data points is 1, which makes it easier to set a suitable regularization parameter $\lambda$. We tried three settings of $\lambda$, as $1/n, 10^{-1}/n$ and $10^{-2}/n$ to represent different levels of regularization. It is expected that when $\lambda = 1/n$, the SVRG/SAGA/SDCA algorithms should converge very fast since the condition number is of the same order as sample size $n$, while for weak regularization case $\lambda = 10^{-2}/n$ these algorithms are expected to converge slowly.

We compare with SVRG (Johnson and Zhang, 2013), SDCA (Shalev-Shwartz and Zhang, 2013), and SAGA Defazio et al. (2014), as they represent popular variance reduced optimization algorithms. we also compare with a related minibatch accelerated SVRG (MB-SVRG) (Nitanda, 2014), which allows large minibatch size without slowing down convergence. For quasi-Newton methods, we compare with L-BFGS with 10 as the limited memory size, and with line search of stepsize to satisfy the Wolfe condition (Wright and Nocedal, 1999). For the proposed MB-SVRP method, at every iteration we simply run one pass of SVRG method initialized with $y_{t-1}$, to approximately solve (4). The implementation is summarized in Algorithm 4, where we fixed many parameters except the minibatch size $b$ as the main tuning parameter. For SVRG and MB-SVRP method, we set $m = 2n/b$ as suggested in (Johnson and Zhang, 2013). Other parameters, such as stepsize in

[2] https://www.csie.ntu.edu.tw/~cjlin/libsvmtools/datasets/

**Algorithm 4** Practical implementation of MB-SVRP Algorithm used in Experiments.

**Parameters**: $b$.
**Initialize** $\eta = \frac{1}{L}$, $\nu = \frac{1-\sqrt{\lambda\eta}}{1+\sqrt{\lambda\eta}}$, $m = \lceil \frac{n}{b} \rceil$, $\widetilde{\lambda} = \frac{1}{\sqrt{b}}$, $\widetilde{w}_0 = 0$.
**Sampling** Sampling $b$ items from $[n]$ to form a minibatch $\bar{B}$.
**for** $s = 1, 2, \ldots$ **do**
  **Calculate** $\widetilde{v} = \frac{1}{n}\sum_{i=1}^{n} \nabla f_i(\widetilde{w}_{s-1})$.
  **Initialize** $y_0 = w_0 = \widetilde{w}_s$.
  **for** $t = 1, 2, \ldots, m$ **do**
    **Sampling** $b$ items from $[n]$ to form a minibatch $B_t$.
    **Initialize** $w_t = y_{t-1}$.
    **Calculate** $\widetilde{u} = \frac{\eta}{b}\sum_{i\in B_t} \nabla f_i(y_{t-1}) + \eta\widetilde{v} - \frac{\eta}{b}\sum_{i\in B_t} \nabla f_i(\widetilde{w}_{s-1})$.
    **for** $k = 1, 2, \ldots, b$ **do**
      **Sampling** Sampling an index $i_k$ from $\bar{B}$.
      **Update** $w_t \leftarrow w_t - \eta\left(\nabla f_{i_k}(w_t) - \nabla f_{i_k}(y_{t-1}) + \widetilde{\lambda}(w_t - y_{t-1}) + \widetilde{u}\right)$.
    **end for**
    **Update**: $y_t = w_t + \nu(w_t - w_{t-1})$.
  **end for**
  **Update** $\widetilde{w}_s = w_m$.
**end for**
**Return** $\widetilde{w}_s$

SVRG and SAGA, and the minibatch size of MP-SVRG, are tuned to give the fastest convergence, the minibatch size of MB-SVRG is set to be the same as MB-SVRP to demonstrate the direct comparison.

Figure 2 and 3 showed the results for $\ell_2$ regularized logistic regression and ridge regression, respectively, where we plot how the objective suboptimality $f(w_t) - f(w^*)$ decreases as the number of gradient evaluations divided by sample size (a.k.a. number of effective passes) increases. We have the following observations:

- When the regularization parameter $\lambda$ (thus strong convexity) is large enough, all methods (except L-BFGS) converges very fast, typically using less than 50 passes over data to converges to numerical precision.
- The convergence of L-BFGS is typically the slowest. When $\lambda$ is large, the advantages of variance-reduced stochastic methods over L-BFGS are significant.
- When the regularization parameter is small, the proposed MP-SVRP method starts to show advantages. In particular, when $\lambda = 10^{-2}/n$, MP-SVRP is substantially faster than all other methods compared in our study.

### 5.2 Empirical Results for Composite Optimization

We also consider empirical comparisons in the setting of non-smooth composition optimization problems. We adopted the same datasets used for the smooth optimization problems, but considered the elastic-net regularized logistic regression and linear regression models (Zou and Hastie, 2005). Where the objective is in the form of (12) but with a elastic-net regularization:

$$g(w) = \frac{\lambda}{2}\|w\|^2 + \mu\|w\|_1,$$

where $\mu$ is set to be $10^{-1}/n$. We follow the same settings as in the previous section, and compare with prox-SVRG (Xiao and Zhang, 2014), prox-SAGA (Defazio et al., 2014) and prox MB-SVRG (Nitanda, 2014), but didn't compare with L-BFGS because the method is generally only applicable for smooth problems. The results are shown in Figure 4 and 5, and we have similar observations as those of the smooth optimization case: all variance-reduced methods converges reasonably fast when $\lambda$ is large, but when $\lambda$ becomes small, the proposed prox MP-SVRP method converges significantly faster because it effectively leveraged the higher-order information from minibatches.

## Acknowledgment

This work was performed when Jialei Wang was a research intern with Tong Zhang at Tencent AI Lab. The authors would like to thank Ji Liu for his helpful feedbacks.

## 6 Conclusion

In this paper, we proposed a novel minibatch stochastic optimization approach for regularized loss minimization in machine learning. This approach efficiently utilizes both variance-reduced first-order gradients and sub-sampled higher order information, and under suitable conditions, it can provably improve the iteration complexity over that of the previous state-of-the-art. Empirical experiments demonstrated that this approach performs well on a variety of smooth and composite optimization tasks in practice. The minibatch nature of the algorithm makes it useful for parallel and distributed computing environment, where additional speed up can be obtained by using many CPU cores for computation.

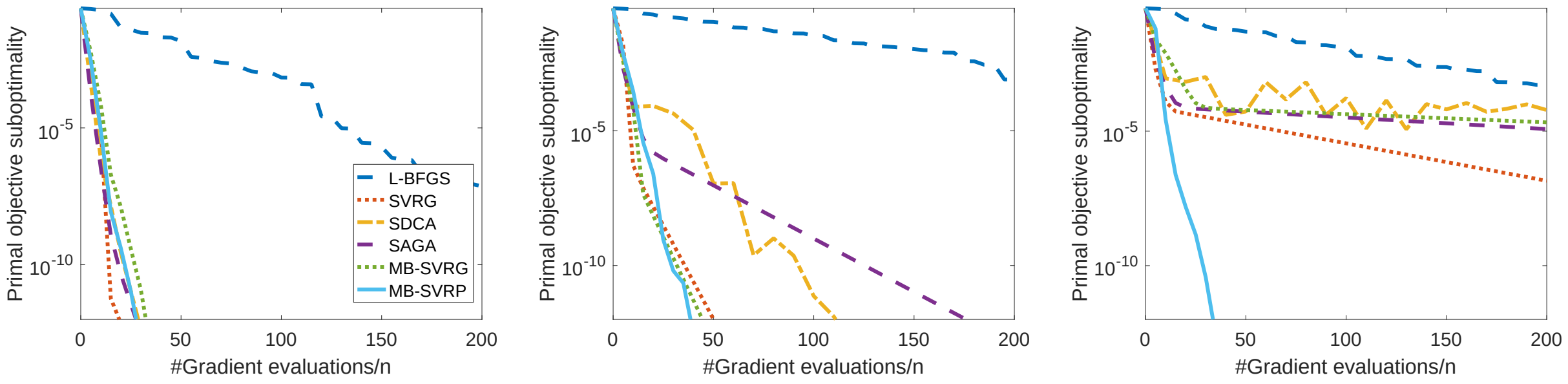

codrna, $\lambda = 1/n$ codrna, $\lambda = 10^{-1}/n$ codrna, $\lambda = 10^{-2}/n$

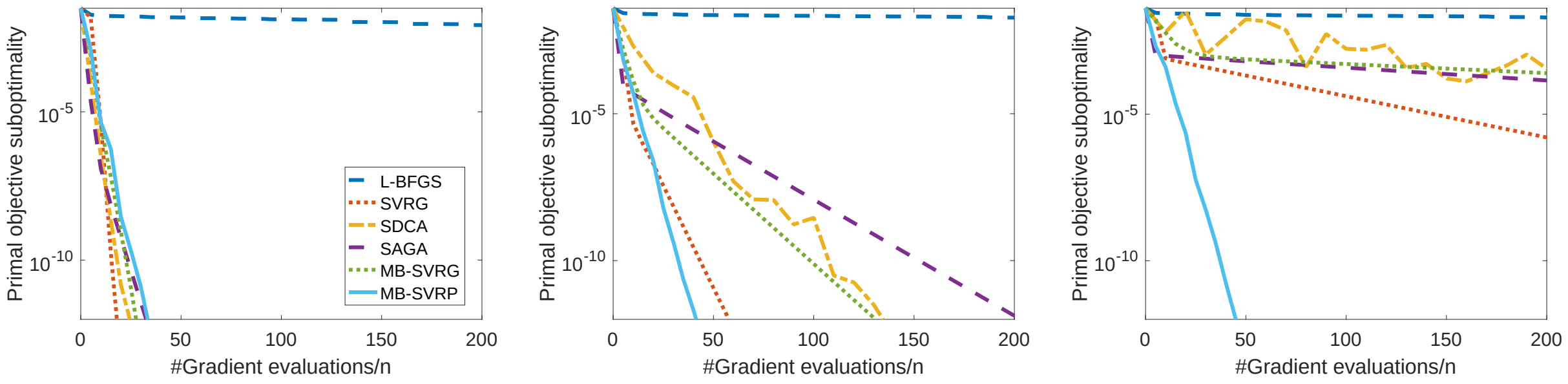

covtype, $\lambda = 1/n$ covtype, $\lambda = 10^{-1}/n$ covtype, $\lambda = 10^{-2}/n$

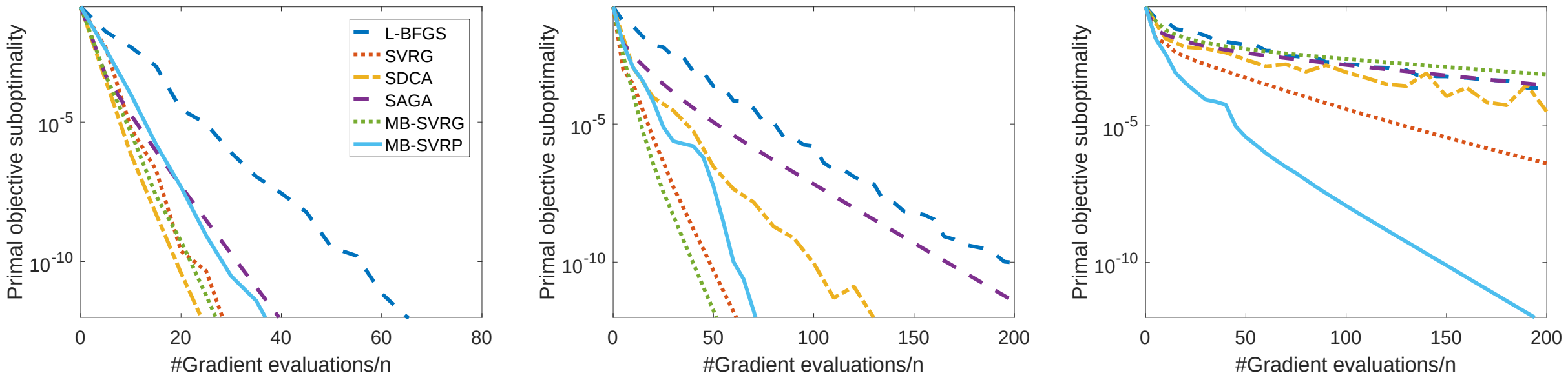

svmguide3, $\lambda = 1/n$ svmguide3, $\lambda = 10^{-1}/n$ svmguide3, $\lambda = 10^{-2}/n$

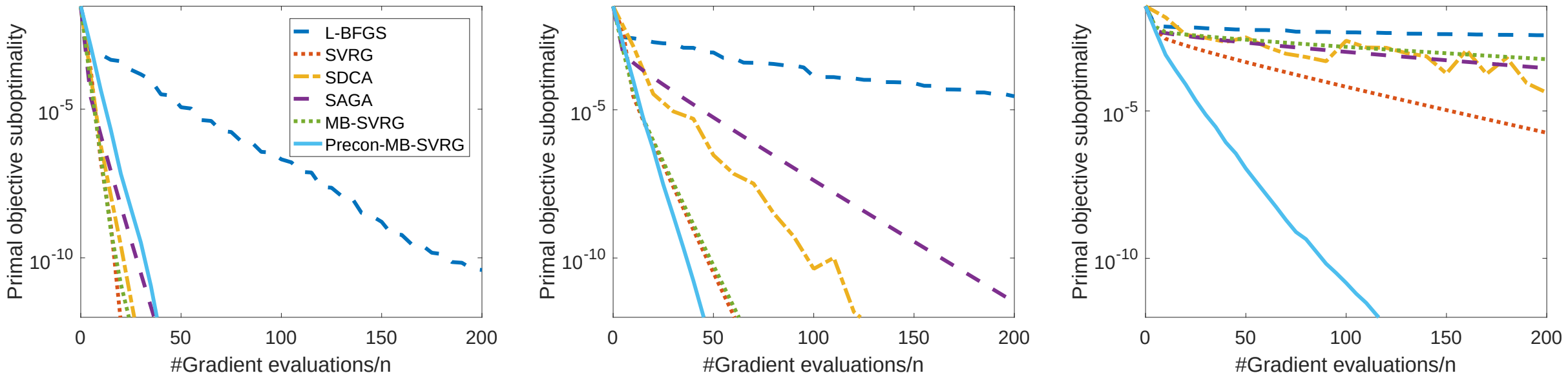

synthetic-c, $\lambda = 1/n$ synthetic-c, $\lambda = 10^{-1}/n$ synthetic-c, $\lambda = 10^{-2}/n$

Figure 2: Comparison of various optimization algorithms for solving $\ell_2$ regularized logistic regression problems.

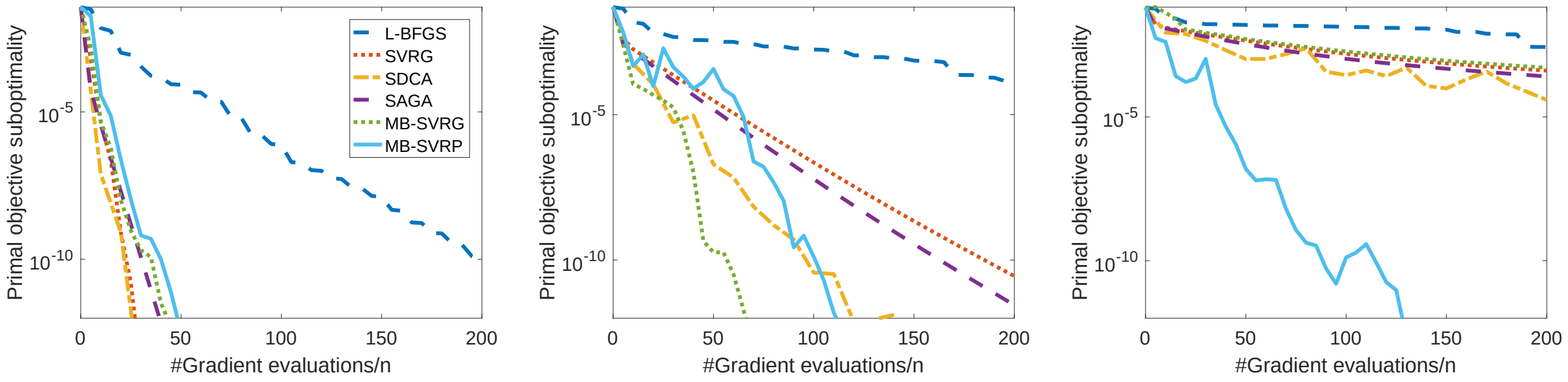

cadata, $\lambda = 1/n$ cadata, $\lambda = 10^{-1}/n$ cadata, $\lambda = 10^{-2}/n$

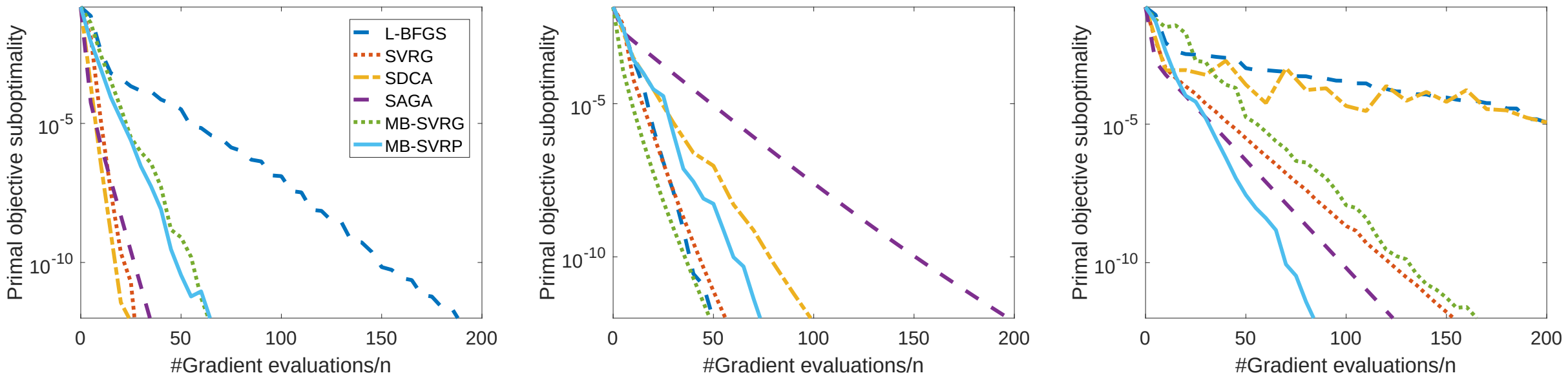

spacega, $\lambda = 1/n$ spacega, $\lambda = 10^{-1}/n$ spacega, $\lambda = 10^{-2}/n$

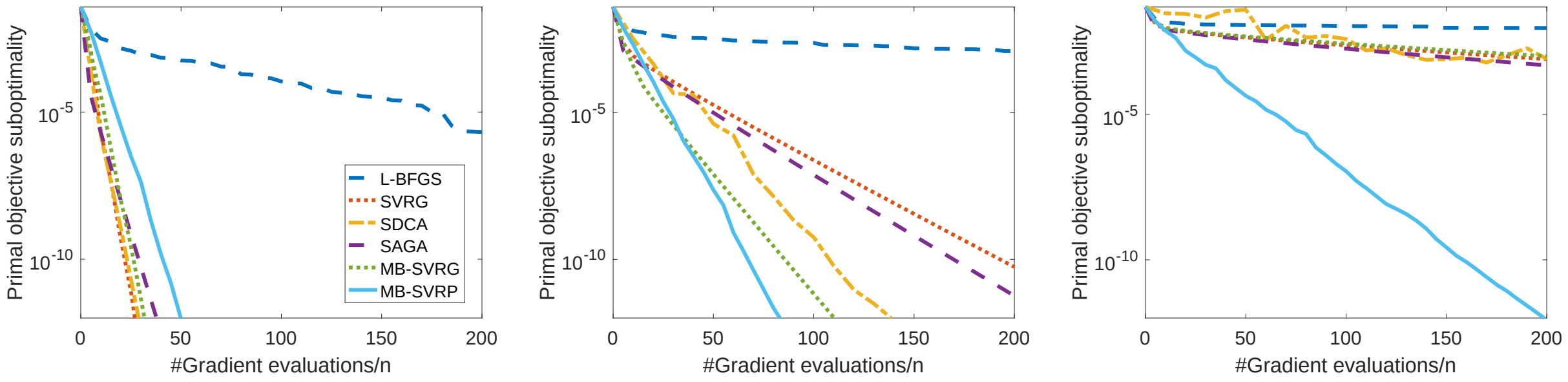

synthetic-r, $\lambda = 1/n$ synthetic-r, $\lambda = 10^{-1}/n$ synthetic-r, $\lambda = 10^{-2}/n$

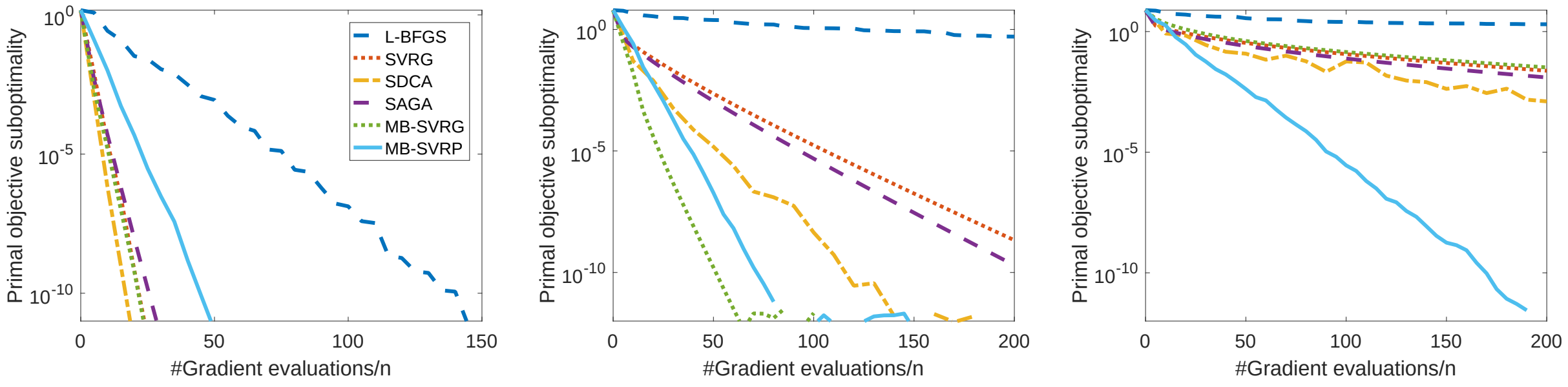

year, $\lambda = 1/n$ year, $\lambda = 10^{-1}/n$ year, $\lambda = 10^{-2}/n$

Figure 3: Comparison of various optimization algorithms for solving ridge regression problems.

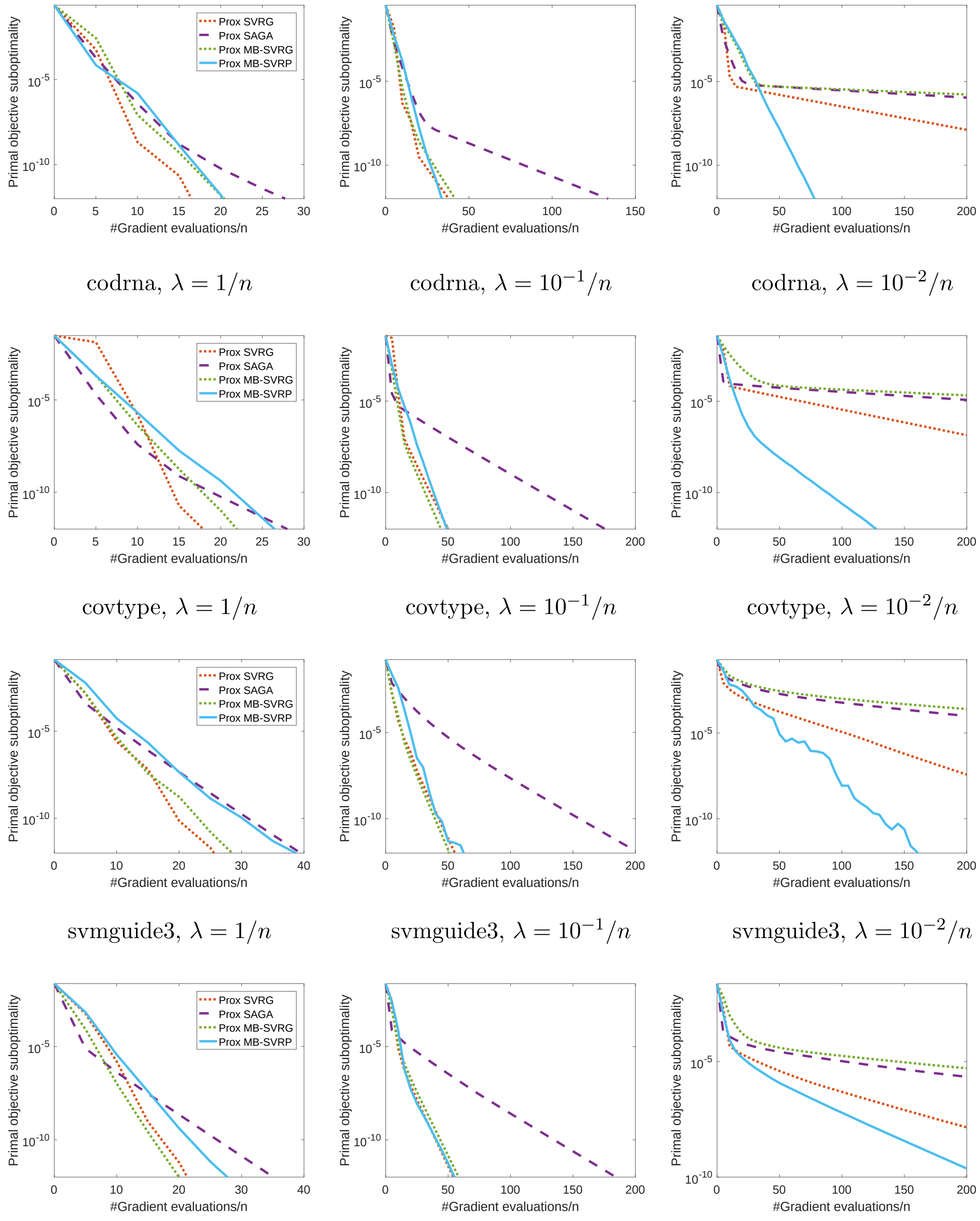

codrna, $\lambda = 1/n$ codrna, $\lambda = 10^{-1}/n$ codrna, $\lambda = 10^{-2}/n$

covtype, $\lambda = 1/n$ covtype, $\lambda = 10^{-1}/n$ covtype, $\lambda = 10^{-2}/n$

svmguide3, $\lambda = 1/n$ svmguide3, $\lambda = 10^{-1}/n$ svmguide3, $\lambda = 10^{-2}/n$

synthetic-c, $\lambda = 1/n$ synthetic-c, $\lambda = 10^{-1}/n$ synthetic-c, $\lambda = 10^{-2}/n$

Figure 4: Comparison of various optimization algorithms for solving elastic-net regularized logistic regression problems.

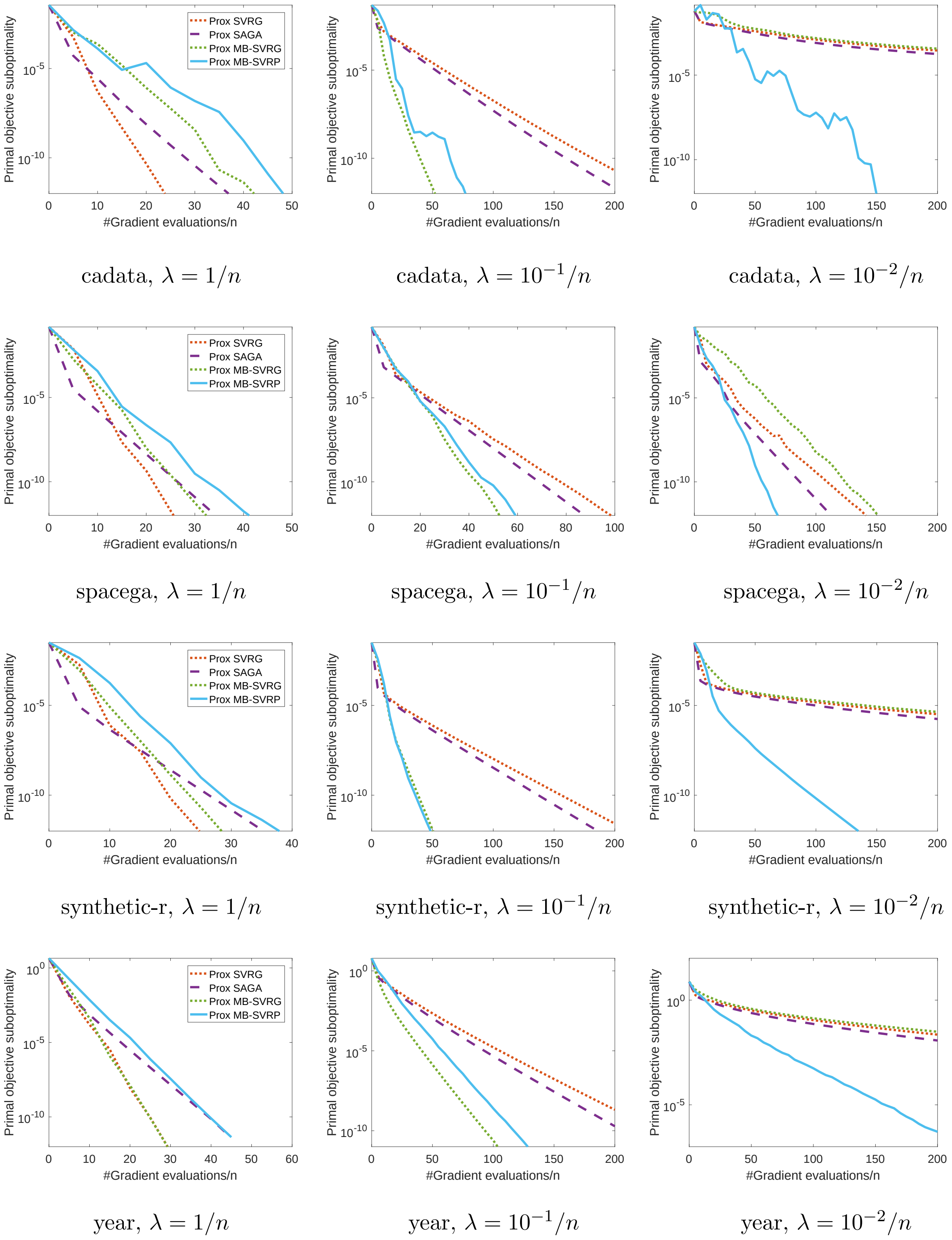

Figure 5: Comparison of various optimization algorithms for solving elastic-net regularized regression problems.

# Appendix

The appendix contains proofs of some theorems and lemmas stated in the main paper.

## A Proofs for Section 3

### A.1 Proof of Proposition 3

*Proof.* Let $w^* = \arg\min_w f(w)$, and $z^* = \arg\min_z f(\widetilde{H}_{\widetilde{\lambda}}^{-1/2}(w^*)z)$, it is clear that $w^* = \widetilde{H}_{\widetilde{\lambda}}^{-1/2}(w^*)z^*$, which is coincident with the exchange of variable rule $z = \widetilde{H}_{\widetilde{\lambda}}^{1/2}(w^*)w$. Moreover, since

$$\begin{aligned}
&\eta\left(\frac{1}{b}\sum_{i\in B_t}\nabla_z f_i(\widetilde{H}_{\widetilde{\lambda}}^{-1/2}(w^*)z_{t-1}) + \frac{1}{n}\sum_{i=1}^{n}\nabla_z f_i(\widetilde{H}_{\widetilde{\lambda}}^{-1/2}(w^*)\widetilde{z}) - \frac{1}{b}\sum_{i\in B_t}\nabla_z f_i(\widetilde{H}_{\widetilde{\lambda}}^{-1/2}(w^*)\widetilde{z})\right)\\
&\overset{①}{=}\eta H_{\widetilde{\lambda}}^{-1/2}(w^*)\left(\frac{1}{b}\sum_{i\in B_t}\nabla f_i(\widetilde{H}_{\widetilde{\lambda}}^{-1/2}(w^*)z_{t-1}) + \frac{1}{n}\sum_{i=1}^{n}\nabla f_i(\widetilde{H}_{\widetilde{\lambda}}^{-1/2}(w^*)\widetilde{z}) - \frac{1}{b}\sum_{i\in B_t}\nabla f_i(\widetilde{H}_{\widetilde{\lambda}}^{-1/2}(w^*)\widetilde{z})\right)\\
&\overset{②}{=}\eta H_{\widetilde{\lambda}}^{-1/2}(w^*)\left(\frac{1}{b}\sum_{i\in B_t}\nabla f_i(w_{t-1}) + \frac{1}{n}\sum_{i=1}^{n}\nabla f_i(\widetilde{w}) - \frac{1}{b}\sum_{i\in B_t}\nabla f_i(\widetilde{w})\right),
\end{aligned}$$

where at step ① we have used the gradient chain rule, at step ② we used the exchange of variables $z = \widetilde{H}_{\widetilde{\lambda}}^{1/2}(w^*)w$. Multiplying both sides of above equation with $\widetilde{H}_{\widetilde{\lambda}}^{-1/2}(w^*)$, and substitute to (15), we recover the exact formulation of (14). □

### A.2 Proof of Lemma 4

*Proof.* We first perform the following decomposition:

$$\|z_t - \bar{z}_t\| \le \|z_t - z_t^*\| + \|z_t^* - \bar{z}_t\|, \tag{22}$$

then we can bound $\|z_t - z_t^*\|$ by

$$\begin{aligned}
\|z_t - z_t^*\| &= \left\|\widetilde{H}_{\widetilde{\lambda}}^{1/2}(w^*)(w_t - w_t^*)\right\| \le \left\|\widetilde{H}_{\widetilde{\lambda}}(w^*)\right\|^{1/2}\|w_t - w_t^*\| \le \sqrt{L+\widetilde{\lambda}}\,\|w_t - w_t^*\|\\
&\le \sqrt{\frac{2(\widetilde{f}_t(w_t) - \widetilde{f}_t(w_t^*))(L+\widetilde{\lambda})}{\lambda+\widetilde{\lambda}}},
\end{aligned}\tag{23}$$

where in the last inequality we use the $\lambda+\widetilde{\lambda}$-strong convexity of $\widetilde{f}_t(w)$. For the term $\|z_t^* - \bar{z}_t\|$, we can bounded it using Hessian Lipschitz property,

$$\begin{aligned}
\|z_t^* - \bar{z}_t\| &= \left\|\eta\left(\widetilde{H}_{\widetilde{\lambda}}(w^*)^{-1} - \widetilde{H}_{\widetilde{\lambda}}(w_{t-1})^{-1}\right)\vartheta_t\right\|\\
&\le \eta\left\|\left(\widetilde{H}_{\widetilde{\lambda}}(w^*)^{-1}\widetilde{H}_{\widetilde{\lambda}}(w_{t-1})^{-1}\right)\left(\widetilde{H}_{\widetilde{\lambda}}(w^*) - \widetilde{H}_{\widetilde{\lambda}}(w_{t-1})\right)\right\|\|\vartheta_t\|\\
&\le \frac{\eta}{\widetilde{\lambda}^2}M\|w_{t-1} - w^*\|\|\vartheta_t\|.
\end{aligned}\tag{24}$$

Combining (22), (23) and (24) finishes the proof. □

### A.3 Proof of Lemma 7

*Proof.* For $j = 1, ..., b$, define random matrix

$$v_j = \frac{1}{b}\left(H_{\widetilde{\lambda}}^{-1/2}(w^*)\ell''(w^{*\top}x_k;b_k)x_kx_k^\top H_{\widetilde{\lambda}}^{-1/2}(w^*) - \frac{1}{n}\sum_{i=1}^{n} H_{\widetilde{\lambda}}^{-1/2}(w^*)\ell''(w^{*\top}x_i;b_i)x_ix_i^\top H_{\widetilde{\lambda}}^{-1/2}(w^*)\right),$$

with probability $1/n, \forall k \in [n]$. It is easy to check that $\mathbb{E}\left[v_j\right] = 0$, and

$$\begin{aligned}
\|v_j\| \leq & \frac{1}{b}\left(\left\|H_{\widetilde{\lambda}}^{-1/2}(w^*)\ell''(w^{*\top}x_k;b_k)x_kx_k^\top H_{\widetilde{\lambda}}^{-1/2}(w^*)\right\|\right) \\
& + \frac{1}{b}\left(\left\|\frac{1}{n}\sum_{i=1}^{n} H_{\widetilde{\lambda}}^{-1/2}(w^*)\ell''(w^{*\top}x_i;b_i)x_ix_i^\top H_{\widetilde{\lambda}}^{-1/2}(w^*)\right\|\right) \\
\leq & \frac{1}{b}\left(2\ell''(w^{*\top}x_i;b_i)\max_{i\in[n]} x_i^\top H_{\widetilde{\lambda}}^{-1}(w^*)x_i\right) \leq \frac{2\rho_{\widetilde{\lambda}}^2}{b}\left(\frac{1}{n}\sum_{i=1}^{n}\ell''(w^{*\top}x_i;b_i)x_i^\top H_{\widetilde{\lambda}}^{-1}(w^*)x_i\right) \\
= & \frac{2\rho_{\widetilde{\lambda}}^2}{b}\operatorname{tr}\left(H_{\widetilde{\lambda}}^{-1}(w^*)\left(\frac{1}{n}\sum_{i=1}^{n}\ell''(w^{*\top}x_i;b_i)x_ix_i^\top\right)\right) = \frac{2\rho_{\widetilde{\lambda}}^2 d_{\widetilde{\lambda}}}{b}.
\end{aligned}$$

for the expected second order moment, denote $\widetilde{x}_k = (\ell''(w^{*\top}x_k;b_k))^{1/2}H_{\widetilde{\lambda}}^{-1/2}(w^*)x_k$ we have

$$\begin{aligned}
\left\|\mathbb{E}\left[v_j^2\right]\right\| = & \left\|\frac{1}{nb^2}\sum_{k=1}^{n}\left(\widetilde{x}_k\widetilde{x}_k^\top - \frac{1}{n}\sum_{i=1}^{n}\widetilde{x}_i\widetilde{x}_i^\top\right)^2\right\| \\
= & \left\|\frac{1}{nb^2}\sum_{k=1}^{n}\left(\left(\widetilde{x}_k\widetilde{x}_k^\top\right)^2 - 2\left(\widetilde{x}_k\widetilde{x}_k^\top\right)\left(\frac{1}{n}\sum_{i=1}^{n}\widetilde{x}_i\widetilde{x}_i^\top\right) + \left(\frac{1}{n}\sum_{i=1}^{n}\widetilde{x}_i\widetilde{x}_i^\top\right)^2\right)\right\| \\
= & \left\|\frac{1}{nb^2}\sum_{k=1}^{n}\left(\widetilde{x}_k\widetilde{x}_k^\top\right)^2 - \frac{2}{b^2}\left(\frac{1}{n}\sum_{i=1}^{n}\widetilde{x}_i\widetilde{x}_i^\top\right)^2 + \frac{1}{b^2}\left(\frac{1}{n}\sum_{i=1}^{n}\widetilde{x}_i\widetilde{x}_i^\top\right)^2\right\| \\
= & \left\|\frac{1}{nb^2}\sum_{k=1}^{n}\left(\widetilde{x}_k\widetilde{x}_k^\top\right)^2 - \frac{1}{b^2}\left(\frac{1}{n}\sum_{i=1}^{n}\widetilde{x}_i\widetilde{x}_i^\top\right)^2\right\| \\
\leq & \left\|\frac{1}{nb^2}\sum_{k=1}^{n}\left(\widetilde{x}_k\widetilde{x}_k^\top\right)^2\right\| \leq \left\|\frac{\rho_{\widetilde{\lambda}}^2 d_{\widetilde{\lambda}}}{b^2}\left(\frac{1}{n}\sum_{k=1}^{n}\widetilde{x}_k\widetilde{x}_k^\top\right)\right\| \\
= & \frac{\rho_{\widetilde{\lambda}}^2 d_{\widetilde{\lambda}}}{b^2}\left\|H_{\widetilde{\lambda}}^{-1/2}(w^*)\left(\frac{1}{n}\sum_{i=1}^{n}\ell''(w^{*\top}x_i;b_i)x_ix_i^\top\right)H_{\widetilde{\lambda}}^{-1/2}(w^*)\right\| \\
\leq & \frac{\rho_{\widetilde{\lambda}}^2 d_{\widetilde{\lambda}}}{b^2}.
\end{aligned}$$

Thus

$$\left\|\sum_{j=1}^{b}\mathbb{E}\left[v_j^2\right]\right\| \leq \frac{\rho_{\widetilde{\lambda}}^2 d_{\widetilde{\lambda}}}{b},$$

using Lemma 23, we know

$$P\left(\left\|H_{\widetilde{\lambda}}^{-1/2}(w^*)\left(\frac{1}{b}\sum_{i\in\bar{B}}\ell''(w^{*\top}x_i;b_i)x_ix_i^\top-\frac{1}{n}\sum_{i=1}^{n}\ell''(w^{*\top}x_i;b_i)x_ix_i^\top\right)H_{\widetilde{\lambda}}^{-1/2}(w^*)\right\|_2\geq t\right)$$

$$=P\left(\left\|\sum_{j=1}^{b}v_j\right\|\geq t\right)\leq d\exp\left(\frac{-t^2/2}{\rho_{\widetilde{\lambda}}^2 d_{\widetilde{\lambda}}/b+2\rho_{\widetilde{\lambda}}^2 d_{\widetilde{\lambda}}t/(3b)}\right),$$

setting $t=2\log\left(\frac{d}{\delta}\right)\cdot\sqrt{\frac{\rho_{\widetilde{\lambda}}^2 d_{\widetilde{\lambda}}}{b}}$ we conclude the proof. □

### A.4 Proof of Lemma 8

*Proof.* When $\widetilde{\lambda}\geq\lambda$, we know

$$\begin{aligned}&\lambda_{\min}(\widetilde{H}_{\widetilde{\lambda}}^{-1/2}(w^*)H_\lambda(w)\widetilde{H}_{\widetilde{\lambda}}^{-1/2}(w^*))\\ &\geq\lambda_{\min}\left(\widetilde{H}_{\widetilde{\lambda}}^{-1/2}(w^*)\left(\frac{\lambda}{\widetilde{\lambda}}\left(\frac{1}{n}\sum_{i=1}^{n}\ell''(w^\top x_i;b_i)x_ix_i^\top\right)+\lambda I\right)\widetilde{H}_{\widetilde{\lambda}}^{-1/2}(w^*)\right)\\ &=\frac{\lambda}{\widetilde{\lambda}}\lambda_{\min}(\widetilde{H}_{\widetilde{\lambda}}^{-1/2}(w^*)H_{\widetilde{\lambda}}(w)\widetilde{H}_{\widetilde{\lambda}}^{-1/2}(w^*)).\end{aligned}$$

Further more we can lower bound $\lambda_{\min}(\widetilde{H}_{\widetilde{\lambda}}^{-1/2}(w^*)H_{\widetilde{\lambda}}(w)\widetilde{H}_{\widetilde{\lambda}}^{-1/2}(w^*))$ by

$$\begin{aligned}&\lambda_{\min}(\widetilde{H}_{\widetilde{\lambda}}^{-1/2}(w^*)H_{\widetilde{\lambda}}(w)\widetilde{H}_{\widetilde{\lambda}}^{-1/2}(w^*))\\ &=\lambda_{\min}(\widetilde{H}_{\widetilde{\lambda}}^{-1}(w^*)H_{\widetilde{\lambda}}(w))\\ &=\frac{1}{\lambda_{\max}(H_{\widetilde{\lambda}}^{-1}(w)\widetilde{H}_{\widetilde{\lambda}}(w^*))}\\ &\geq\frac{1}{\lambda_{\max}(H_{\widetilde{\lambda}}^{-1}(w)H_{\widetilde{\lambda}}(w))+\lambda_{\max}(H_{\widetilde{\lambda}}^{-1}(w)(H_{\widetilde{\lambda}}(w)-\widetilde{H}_{\widetilde{\lambda}}(w^*)))}\\ &\geq\frac{1}{1+\lambda_{\max}(H_{\widetilde{\lambda}}^{-1}(w)(H_{\widetilde{\lambda}}(w^*)-\widetilde{H}_{\widetilde{\lambda}}(w^*)))+\lambda_{\max}(H_{\widetilde{\lambda}}^{-1}(w)(H_{\widetilde{\lambda}}(w)-H_{\widetilde{\lambda}}(w^*)))},\end{aligned}\tag{25}$$

For the term $\lambda_{\max}(H_{\widetilde{\lambda}}^{-1}(w)(H_{\widetilde{\lambda}}(w^*)-\widetilde{H}_{\widetilde{\lambda}}(w^*)))$, we have the following decomposition:

$$\begin{aligned}\lambda_{\max}(H_{\widetilde{\lambda}}^{-1}(w)(H_{\widetilde{\lambda}}(w^*)-\widetilde{H}_{\widetilde{\lambda}}(w^*)))\leq&\lambda_{\max}(H_{\widetilde{\lambda}}^{-1}(w^*)(H_{\widetilde{\lambda}}(w^*)-\widetilde{H}_{\widetilde{\lambda}}(w^*)))\\ &+\lambda_{\max}((H_{\widetilde{\lambda}}^{-1}(w^*)-H_{\widetilde{\lambda}}^{-1}(w))(H_{\widetilde{\lambda}}(w^*)-\widetilde{H}_{\widetilde{\lambda}}(w^*)))\end{aligned}\tag{26}$$

since with probability at least $1-\delta$, we have

$$\lambda_{\max}(H_{\widetilde{\lambda}}^{-1}(w^*)(H_{\widetilde{\lambda}}(w^*)-\widetilde{H}_{\widetilde{\lambda}}(w^*)))\leq 2\log\left(\frac{d}{\delta}\right)\sqrt{\frac{\rho_{\widetilde{\lambda}}^2 d_{\widetilde{\lambda}}}{b}},\tag{27}$$

and

$$\lambda_{\max}((H_{\widetilde{\lambda}}^{-1}(w^*)-H_{\widetilde{\lambda}}^{-1}(w))(H_{\widetilde{\lambda}}(w^*)-\widetilde{H}_{\widetilde{\lambda}}(w^*)))\leq\lambda_{\max}(H_{\widetilde{\lambda}}^{-1}(w^*)-H_{\widetilde{\lambda}}^{-1}(w))L\leq ML\left\|w-w^*\right\|\tag{28}$$

by Lemma (7), and by Lipschitz Hessian property (Condition 1). We also have

$$\lambda_{\max}(H_{\widetilde{\lambda}}^{-1}(w)(H_{\widetilde{\lambda}}(w) - H_{\widetilde{\lambda}}(w^*))) \leq \frac{\lambda_{\max}(H_{\widetilde{\lambda}}(w) - H_{\widetilde{\lambda}}(w^*))}{\widetilde{\lambda}} \leq \frac{M\left\|w - w^*\right\|}{\widetilde{\lambda}}. \tag{29}$$

Combining (25), (27), (28), (29) we have with probability at least $1-\delta$

$$\min_{\|w-w^*\|\leq R} \lambda_{\min}(\widetilde{H}_{\widetilde{\lambda}}^{-1/2}(w^*)H_{\lambda}(w)\widetilde{H}_{\widetilde{\lambda}}^{-1/2}(w^*)) \geq \frac{\lambda}{\widetilde{\lambda}} \cdot \frac{1}{1 + 2\log\left(\frac{d}{\delta}\right)\sqrt{\frac{\rho_{\widetilde{\lambda}}^2 d_{\widetilde{\lambda}}}{b}} + MLR + \frac{MR}{\widetilde{\lambda}}},$$

which finishes the proof. □

## A.5 Proof of Lemma 9

*Proof.* First we upper bound $\lambda_{\max}\left(\widetilde{H}_{\widetilde{\lambda}}^{-1/2}(w^*)H_{\widetilde{\lambda}}(w^*)\widetilde{H}_{\widetilde{\lambda}}^{-1/2}(w^*)\right)$, since when $b \geq 16\rho_{\widetilde{\lambda}}^2 d_{\widetilde{\lambda}} \log^2(d/\delta)$ we have $\lambda_{\max}\left(H_{\widetilde{\lambda}}^{-1}(w^*)(H_{\widetilde{\lambda}}(w^*) - \widetilde{H}_{\widetilde{\lambda}}(w^*))\right) \leq 1$, thus

$$\begin{aligned}
\lambda_{\max}\left(\widetilde{H}_{\widetilde{\lambda}}^{-1/2}(w^*)H_{\widetilde{\lambda}}(w^*)\widetilde{H}_{\widetilde{\lambda}}^{-1/2}(w^*)\right) =& \frac{1}{\lambda_{\min}\left(H_{\widetilde{\lambda}}^{-1}(w^*)\widetilde{H}_{\widetilde{\lambda}}(w^*)\right)} \\
\leq& \frac{1}{1 - \lambda_{\max}\left(H_{\widetilde{\lambda}}^{-1}(w^*)(H_{\widetilde{\lambda}}(w^*) - \widetilde{H}_{\widetilde{\lambda}}(w^*))\right)} \\
\leq& \frac{1}{1 - 2\log\left(\frac{d}{\delta}\right)\sqrt{\frac{\rho_{\widetilde{\lambda}}^2 d_{\widetilde{\lambda}}}{b}}} \\
\leq& 2.
\end{aligned}$$

Then we bound $\max_{i\in[n]} \ell''(w^{*\top}x_i; b_i)x_i^\top \widetilde{H}_{\widetilde{\lambda}}^{-1}(w^*)x_i$ through $\max_{i\in[n]} \ell''(w^{*\top}x_i; b_i)x_i^\top H_{\widetilde{\lambda}}^{-1}(w^*)x_i$, because

$$\begin{aligned}
x_i^\top \widetilde{H}_{\widetilde{\lambda}}^{-1}(w^*)x_i =& x_i^\top H_{\widetilde{\lambda}}^{-1}(w^*)x_i + x_i^\top(\widetilde{H}_{\widetilde{\lambda}}^{-1}(w^*) - H_{\widetilde{\lambda}}^{-1}(w^*))x_i \\
=& x_i^\top H_{\widetilde{\lambda}}^{-1}(w^*)x_i + x_i^\top(\widetilde{H}_{\widetilde{\lambda}}^{-1}(w^*)H_{\widetilde{\lambda}}(w^*) - I)H_{\widetilde{\lambda}}^{-1}(w^*)x_i \\
\leq& x_i^\top H_{\widetilde{\lambda}}^{-1}(w^*)x_i + \lambda_{\max}(\widetilde{H}_{\widetilde{\lambda}}^{-1}(w^*)H_{\widetilde{\lambda}}(w^*) - I)x_i^\top H_{\widetilde{\lambda}}^{-1}(w^*)x_i \\
\overset{①}{\leq}& 2x_i^\top H_{\widetilde{\lambda}}^{-1}(w^*)x_i,
\end{aligned}$$

where in step ① we used the fact that

$$\lambda_{\max}(\widetilde{H}_{\widetilde{\lambda}}^{-1}(w^*)H_{\widetilde{\lambda}}(w^*) - I) \leq \max\left\{|\lambda_{\max}(\widetilde{H}_{\widetilde{\lambda}}^{-1}(w^*)H_{\widetilde{\lambda}}(w^*)) - 1|, |\lambda_{\min}(\widetilde{H}_{\widetilde{\lambda}}^{-1}(w^*)H_{\widetilde{\lambda}}(w^*)) - 1|\right\} \leq 1.$$

Based on the definition of effective dimension (Definition 5) and condition of bounded statistical leverage (Assumption 6), we can bound the smoothness as

$$
\begin{aligned}
\max_{i\in[n]} \ell''({w^*}^\top x_i; b_i) x_i^\top \widetilde{H}_{\widetilde{\lambda}}^{-1}(w^*) x_i \leq& 2\max_{i\in[n]} \ell''({w^*}^\top x_i; b_i) x_i^\top H_{\widetilde{\lambda}}^{-1}(w^*) x_i \\
\leq& 2\rho_{\widetilde{\lambda}}^2 \left(\frac{1}{n}\sum_{i=1}^n \ell''({w^*}^\top x_i; b_i) x_i^\top H_{\widetilde{\lambda}}^{-1}(w^*) x_i\right) \\
=& 2\rho_{\widetilde{\lambda}}^2 \left(\frac{1}{n}\sum_{i=1}^n \mathrm{tr}(\ell''({w^*}^\top x_i; b_i) x_i^\top H_{\widetilde{\lambda}}^{-1}(w^*) x_i)\right) \\
\overset{①}{=}& 2\rho_{\widetilde{\lambda}}^2 \left(\frac{1}{n}\sum_{i=1}^n \mathrm{tr}(\ell''({w^*}^\top x_i; b_i) x_i x_i^\top H_{\widetilde{\lambda}}^{-1}(w^*))\right) \\
=& 2\rho_{\widetilde{\lambda}}^2 \left(\mathrm{tr}\left(\frac{1}{n}\sum_{i} \ell''({w^*}^\top x_i; b_i) x_i x_i^\top H_{\widetilde{\lambda}}^{-1}(w^*)\right)\right) \\
=& 2\rho_{\widetilde{\lambda}}^2 \left(\mathrm{tr}\left(\left(\frac{1}{n}\sum_{i} \ell''({w^*}^\top x_i; b_i) x_i x_i^\top\right) H_{\widetilde{\lambda}}^{-1}(w^*)\right)\right) \\
=& 2\rho_{\widetilde{\lambda}}^2 \sum_{j=1}^d \frac{\lambda_j}{\lambda_j + \widetilde{\lambda}} = 2\rho_{\widetilde{\lambda}}^2 d_{\widetilde{\lambda}},
\end{aligned}
$$

which finished the proof, where in step ① we have used the fact that $\mathrm{tr}(ABC) = \mathrm{tr}(CAB)$ for any $A, B, C$. □

### A.6 Proof of Theorem 10

*Proof.* Since we choose $\widetilde{\lambda}$ at the scale of $\widetilde{\lambda} = \max\left\{\lambda, \frac{L}{b}\right\}$, then we can apply Lemma 8 and Lemma 9 to verify that the new condition number after "preconditioning" will be

$$
\max\left\{4\rho_{\widetilde{\lambda}}^2 d_{\widetilde{\lambda}}, \frac{4L\rho_{\widetilde{\lambda}}^2 d_{\widetilde{\lambda}}}{\lambda b}\right\},
$$

applying Corollary 14 we know as long as the inexactness condition (53) is satisfied, we requires $\mathcal{O}\left(\log\left(\frac{1}{\epsilon}\right)\right)$ full gradient evaluations and

$$
\max\left\{C_{\mathrm{sb}}, C_{\mathrm{lb}}\right\} \cdot \frac{10}{9}\log(36)\log\left(\frac{1}{\epsilon}\right) \tag{30}
$$

total calls of approximate minimization of (8) to ensure $\mathbb{E}f(\widetilde{w}_s) - f(w^*) \leq \epsilon$, where the factors $C_{\mathrm{sb}}, C_{\mathrm{lb}}$ are

$$
\begin{aligned}
C_{\mathrm{sb}} =& \mathcal{O}\left(\max\left\{\frac{240}{b}, \frac{240L}{\lambda b^2}\right\} \cdot \rho_{\widetilde{\lambda}}^2 d_{\widetilde{\lambda}}\right), \\
C_{\mathrm{lb}} =& \mathcal{O}\left(\max\left\{2\sqrt{6}, 2\sqrt{6}\sqrt{\frac{L}{\lambda b}}\right\} \cdot \sqrt{\rho_{\widetilde{\lambda}}^2 d_{\widetilde{\lambda}}}\right),
\end{aligned}
$$

which represent the cases of small and large minibatch sizes, respectively. Since the condition number of (8) is never larger than $L \cdot \frac{b}{L} = b$, we know from Lemma 24, when applying SVRG to solve (8), to reach some point of which objective the suboptimality of (8) satisfies

$$w_t - \min_w \widetilde{f}_t(w) \le \frac{1}{10^5} \cdot \left(\frac{\lambda}{L}\right)^7 \epsilon, \tag{31}$$

the following number of gradient calls sufficient:

$$C \cdot b \cdot \log\left(\frac{L}{\lambda}\right) \cdot \log\left(\frac{1}{\epsilon}\right), \tag{32}$$

for some universal constant $C$. On the other hand, for every subproblem (8), if (31) is satisfied, applying Lemma 4 and Lemma 21 we know the gradient error in IMBA-SVRG (Algorithm 3)in the preconditioned space can be upper bounded by

$$\|\xi_t\|^2 \le \left(\frac{\lambda}{L}\right)^6 \cdot \frac{2\epsilon}{10^5} + \frac{\eta^2 M^2 R^2 L}{\widetilde{\lambda}^4}(f(\widetilde{w}_{s-1}) - f(w^*)),$$

thus condition (53) in Theorem 13 is satisfied. Combining (30) and (32) we know the total number of gradient calls:

$$n \cdot \log\left(\frac{1}{\epsilon}\right) + \frac{10C}{9} \log(36) \cdot \max\{C_{\text{sb}}, C_{\text{lb}}\} \cdot b \cdot \log\left(\frac{L}{\lambda}\right) \cdot \log^2\left(\frac{1}{\epsilon}\right) \tag{33}$$

is sufficient to obtain a solution such that $\mathbb{E} f(\widetilde{w}_s) - f(w^*) \le \epsilon$ is satisfied. Next we choose $b$ such that the total iteration complexity of above expression is minimized. Start from here we will ignore the constants before these factors for simplicity, we know the term $\max\{C_{\text{sb}}, C_{\text{lb}}\} \cdot b$ is of order

$$\mathcal{O}\left(\max\left\{\rho_{\widetilde{\lambda}}^2 d_{\widetilde{\lambda}}, \frac{\rho_{\widetilde{\lambda}}^2 d_{\widetilde{\lambda}} L}{\lambda b}, \sqrt{\frac{\rho_{\widetilde{\lambda}}^2 d_{\widetilde{\lambda}} L b}{\lambda}}\right\}\right),$$

when $n \gtrsim \rho_{\widetilde{\lambda}}^2 d_{\widetilde{\lambda}} \left(\frac{L}{\lambda}\right)^{1/3}$, we can choose $b \asymp \rho_{\widetilde{\lambda}}^2 d_{\widetilde{\lambda}} \left(\frac{L}{\lambda}\right)^{1/3}$, then $\max\{C_{\text{sb}}, C_{\text{lb}}\} \cdot b$ is of order

$$\mathcal{O}\left(\max\left\{\rho_{\widetilde{\lambda}}^2 d_{\widetilde{\lambda}}, \left(\frac{L}{\lambda}\right)^{2/3}, \rho_{\widetilde{\lambda}}^2 d_{\widetilde{\lambda}} \left(\frac{L}{\lambda}\right)^{2/3}\right\}\right) = \mathcal{O}\left(\rho_{\widetilde{\lambda}}^2 d_{\widetilde{\lambda}} \left(\frac{L}{\lambda}\right)^{2/3}\right), \tag{34}$$

when $n \lesssim \rho_{\widetilde{\lambda}}^2 d_{\widetilde{\lambda}} \left(\frac{L}{\lambda}\right)^{1/3}$, we can choose $b \asymp n$, in which case $\max\{C_{\text{sb}}, C_{\text{lb}}\} \cdot b$ can be upper bounded by

$$\mathcal{O}\left(\max\left\{\frac{\rho_{\widetilde{\lambda}}^2 d_{\widetilde{\lambda}} L}{\lambda n}, \sqrt{\frac{\rho_{\widetilde{\lambda}}^2 d_{\widetilde{\lambda}} L n}{\lambda}}\right\}\right) = \mathcal{O}\left(\frac{\rho_{\widetilde{\lambda}}^2 d_{\widetilde{\lambda}} L}{\lambda n}\right), \tag{35}$$

Combining (33), (34) and (35), we know the total number of individual function gradient calls to reach $\epsilon$-suboptimality is

$$\mathcal{O}\left(\widetilde{\kappa} \cdot \log\left(\frac{L}{\lambda}\right) \cdot \log^2\left(\frac{1}{\epsilon}\right) + n \cdot \log\left(\frac{1}{\epsilon}\right)\right),$$

where

$$\widetilde{\kappa} = \max\left\{\rho_{\widetilde{\lambda}}^2 d_{\widetilde{\lambda}} \left(\frac{L}{\lambda}\right)^{2/3}, \frac{\rho_{\widetilde{\lambda}}^2 d_{\widetilde{\lambda}} L}{\lambda n}\right\},$$

which finishes the proof. □

### A.7 Proof of Theorem 12

*Proof.* Applying the theory of catalyst acceleration (Lemma 25) we know only $\mathcal{O}\left(\left(\sqrt{\frac{\lambda+\gamma}{\lambda}}\right)\log\left(\frac{1}{\epsilon}\right)\right)$ calls of MP-SVRP is sufficient to reach $\epsilon$-objective suboptimality, as long as each iterate $w_r$ satisfies

$$f(w_r) + \frac{\gamma}{2}\|w_r - z_{r-1}\|^2 \leq \min_w f(w) + \frac{\gamma}{2}\|w - z_{r-1}\|^2 + \frac{\lambda\epsilon}{3600(\lambda+\gamma)}\left(1 - \frac{9}{10}\sqrt{\frac{\lambda}{\lambda+\gamma}}\right).$$

Moreover, according to Theorem 10, the iteration complexity of solving a $\lambda+\gamma$ strongly convex problem (20) is

$$\widetilde{\mathcal{O}}\left(n + \max\left\{\rho_{\widetilde{\lambda}}^2 d_{\widetilde{\lambda}}\left(\frac{L}{(\lambda+\gamma)}\right)^{2/3}, \frac{\rho_{\widetilde{\lambda}}^2 d_{\widetilde{\lambda}} L}{(\lambda+\gamma)n}\right\}\right),$$

combining these two results we know the total iteration complexity of Algorithm 2 can be upper bounded by

$$\widetilde{\mathcal{O}}\left(\left(\sqrt{\frac{\lambda+\gamma}{\lambda}}\right)\left(n + \max\left\{\left(\frac{\rho_{\widetilde{\lambda}}^2 d_{\widetilde{\lambda}} L}{(\lambda+\gamma)}\right)^{2/3}, \frac{\rho_{\widetilde{\lambda}}^2 d_{\widetilde{\lambda}} L}{(\lambda+\gamma)n}\right\}\right)\right),$$

When $\frac{\rho_{\widetilde{\lambda}}^2 d_{\widetilde{\lambda}} L}{\lambda} \geq n^{3/2}$, if we choose

$$\gamma = \frac{\rho_{\widetilde{\lambda}}^2 d_{\widetilde{\lambda}} L}{n^{3/2}} - \lambda,$$

we obtain the iteration complexity can be upper bounded by

$$\widetilde{\mathcal{O}}\left(\left(\sqrt{\frac{\rho_{\widetilde{\lambda}}^2 d_{\widetilde{\lambda}} L}{n^{3/2}\lambda}}\right)\cdot n\right) = \widetilde{\mathcal{O}}\left(\sqrt{\rho_{\widetilde{\lambda}}^2 d_{\widetilde{\lambda}}}\cdot n^{1/4}\left(\frac{L}{\lambda}\right)^{1/2}\right),$$

which concludes the proof. □

## B Proofs for Section 4

In this appendix, we describe the stochastic estimation sequence approach for inexact accelerated minibatch SVRG in Section 4, and provide the detailed analysis.

### B.1 Stochastic Estimation Sequences

We introduce the following definition.

**Definition 15.** ***(Stochastic estimation sequence)*** *A sequence of pairs $\{V_t(w), \theta_t\}_{t\geq 0}$ is called an estimation sequence of the function $f(w)$ if $\theta_t > 0$ and for any $w \in \mathbb{R}^d, t \geq 0$ we have*

$$V_t(w) \leq (1-\theta_t)f(w) + \theta_t V_0(w), \tag{36}$$

*if $\{V_t(w)\}_{t\geq 0}$ is a sequence of random functions, we call the sequence of pairs $\{V_t(w), \theta_t\}_{t\geq 0}$ a stochastic estimation sequence if (36) holds in expectation, i.e.*

$$\mathbb{E}\left[V_t(w)\right] \leq (1-\theta_t)f(w) + \theta_t V_0(w).$$

The following lemma is a generalized version of Lemma 1 in (Lin et al., 2014), which shows if we can construct upper bound of $\mathbb{E}[f(w_t)]$ using $\mathbb{E}[V_t(w)]$, then we get convergence of $\mathbb{E}[f(w_t)]$.

**Lemma 16.** *Suppose $\{V_t(w), \theta_t\}_{t\geq 0}$ is a stochastic estimation sequence of the function $f(w)$. Let $w^*$ be the minimizer of $f(w)$. If there are sequences of random variables $\{w_t\}_{t\geq 0}$ and $\{\epsilon_t\}_{t\geq 0}$ in $\mathbb{R}^d$, $\{\delta_t\}_{t\geq 0}$ in $\mathbb{R}$ such that*

$$\mathbb{E}[f(w_t)] \leq \min_w \{\mathbb{E}[V_t(w)]\} + \mathbb{E}[\delta_t] \tag{37}$$

*holds for all $t \geq 0$, then*

$$\mathbb{E}[f(w_t)] - f(w^*) \leq \theta_t(V_0(w^*) - f(w^*)) + \mathbb{E}[\delta_t].$$

The following lemma constructs a concrete stochastic estimation sequence.

**Lemma 17.** *Assume $f(w) = \frac{1}{n}\sum_{i=1}^n f_i(w)$, where each $f_i(w)$ is $\lambda$-strongly convex and $L$-smooth. Suppose that*

- *$V_0(w)$ is an arbitrary deterministic function on $\mathbb{R}^d$;*
- *$\{\alpha_t\}_{t\geq 0}$ is a sequence that satisfies $\alpha_t \in (0,1), \forall t \geq 0$, and $\sum_{t=0}^\infty \alpha_t = \infty$.*
- *$\{y_t\}_{t=0}^\infty$ is an arbitrary sequence in $\mathbb{R}^d$;*
- *Define $\{v_t\}_{t\geq 1}$ as:*

$$v_t = \frac{1}{b}\sum_{i\in B_t} \nabla f_i(y_{t-1}) - \frac{1}{b}\sum_{i\in B_t} \nabla f_i(\widetilde{w}_{s-1}) + \frac{1}{n}\sum_{i=1}^n \nabla f_i(\widetilde{w}_{s-1}).$$

*Define the sequence $\{w_t\}_{t\geq 0}$ and $\{V_t(w)\}_{t\geq 0}$ as follows. Let $y_0 = w_0$ and $w_t$ be arbitrary vector in $\mathbb{R}^d$ such that*

$$w_t = y_{t-1} - \eta v_t + \xi_t,$$

*where the stepsize $\eta$ satisfies $0 < \eta \leq \frac{1}{L}$. Let*

$$\begin{aligned} V_t(w) =& (1-\alpha_{t-1})V_{t-1}(w) \\ &+ \alpha_{t-1}\left(\frac{1}{b}\sum_{i\in B_t} f_i(y_{t-1}) + \left\langle v_t - \frac{\xi_t}{\eta}, w - y_{t-1}\right\rangle + \frac{\lambda}{4}\|w - y_{t-1}\|^2 - \frac{\|\xi_t\|^2}{\lambda\eta^2}\right), \end{aligned} \tag{38}$$

*let $\theta_0 = 1$, and $\theta_t = (1-\alpha_{t-1})\theta_{t-1}, \forall t \geq 1$. Then the sequence $\{V_t(w), \theta_t\}_{t\geq 0}$ is a stochastic estimation sequence of $f(w)$.*

By above construction, we know $V_t(w)$ is a quadratic function, its minimizer and minimum value has the following iterative form:

**Lemma 18.** *Define* $\forall t \geq 0$*:*

$$V_t^* = \min_w V_t(w),$$

*suppose we choose the function* $V_0(w)$ *in Lemma 17 as*

$$V_0(w) = V_0^* + \frac{\lambda}{4}\left\|w - z_0\right\|^2,$$

*with* $V_0^* = f(z_0), \epsilon_0 = 0$ *and* $\delta_0 = 0$ *for* $z_0 = w_0$*. Then the sequence* $\{V_t(w)\}_{t\geq 0}$ *defined in 38 can be written as*

$$V_t(w) = V_t^* + \frac{\lambda}{4}\left\|w - z_t\right\|^2, \tag{39}$$

*where the sequences* $\{V_t^*\}, \{z_t\}$ *are defined recursively as:*

$$z_t = (1-\alpha_{t-1})z_{t-1} + \alpha_{t-1}y_{t-1} - \frac{2\alpha_{t-1}}{\lambda}\left(v_t - \frac{\xi_t}{\eta}\right), \tag{40}$$

$$\begin{aligned} V_t^* =& (1-\alpha_{t-1})V_{t-1}^* + \frac{\alpha_{t-1}(1-\alpha_{t-1})\lambda}{4}\left\|z_{t-1} - y_{t-1}\right\|^2 + \alpha_{t-1}(1-\alpha_{t-1})\left\langle v_t - \frac{\xi_t}{\eta}, z_{t-1} - y_{t-1}\right\rangle \\ &- \frac{\alpha_{t-1}^2}{\lambda}\left\|v_t\right\|^2 + \frac{\alpha_{t-1}}{b}\sum_{i\in B_t} f_i(y_{t-1}) - \frac{(\alpha_{t-1} + \alpha_{t-1}^2)\left\|\xi_t\right\|^2}{\lambda\eta^2} + \frac{2\alpha_{t-1}^2}{\lambda\eta}\langle \xi_t, v_t\rangle. \end{aligned} \tag{41}$$

The following lemma establishes a connection between $z_t, y_t, w_t$ when $z_t$ is updated as (40), and $y_t, w_t$ is updated as Algorithm 3.

**Lemma 19.** *Based on the updating rule of* $z_t$ *in* (40)*, if we choose* $\alpha_t = \alpha = \sqrt{\frac{\lambda\eta}{2}}, \forall t \geq 0$*, we have the following inequalities hold for all* $t \geq 0$*:*

$$z_t - y_t = \frac{1}{\alpha}(y_t - w_t).$$

The following lemma gives concrete expression of $\delta_t$ such that for the IMBA-SVRG algorithm, the condition (37) in Lemma 16 is satisfied.

**Lemma 20.** *Suppose we choose* $\forall t \geq 0$

$$\alpha_t = \alpha = \sqrt{\frac{\lambda\eta}{2}},$$

*and the stepsize* $\eta$ *satisfying* $\eta \leq \frac{1}{8L}$*, we have* $\forall t \geq 0$

$$\mathbb{E}\left[f(w_t)\right] \leq \mathbb{E}\left[V_t^*\right] + \mathbb{E}\left[\delta_t\right],$$

*where* $\delta_0 = 0$ *and* $\forall t \geq 1$*:*

$$\delta_t = \left(1 - \sqrt{\frac{\lambda\eta}{2}}\right)\delta_{t-1} + \left(\eta\left\|v_t - \nabla f(y_{t-1})\right\|^2 - \frac{(1-\sqrt{\lambda\eta/2})\lambda}{\sqrt{2\lambda\eta}}\left\|y_{t-1} - w_{t-1}\right\|^2 + \frac{1}{4\lambda^2\eta^3}\left\|\xi_t\right\|^2\right).$$

The lemma below bounds the term $\mathbb{E}\left[\eta \left\|v_t - \nabla f(y_{t-1})\right\|^2 - \frac{(1-\sqrt{\lambda\eta/2})\lambda}{\sqrt{2\lambda\eta}} \left\|y_{t-1} - w_{t-1}\right\|^2\right]$.

**Lemma 21.** *Suppose $B_t$ is constructed by uniform sampling with or without replacement, we have the following inequality holds:*

$$\begin{aligned}
&\mathbb{E}\left[\eta \left\|v_t - \nabla f(y_{t-1})\right\|^2 - \frac{(1-\sqrt{\lambda\eta/2})\lambda}{\sqrt{2\lambda\eta}} \left\|y_{t-1} - w_{t-1}\right\|^2\right] \\
&\leq \frac{8\eta L}{b}\left(f(w_{t-1}) - f(w^*) + f(\widetilde{w}_{s-1}) - f(w^*)\right) + \left(\frac{2\eta L^2}{b} - \frac{(1-\sqrt{\lambda\eta/2})\lambda}{\sqrt{2\lambda\eta}}\right)\left\|y_{t-1} - w_{t-1}\right\|^2.
\end{aligned}$$

Based on above lemma, by carefully choosing the stepsize $\eta$, we get the following lemma which is important to obtain iteration complexity for IMBA-SVRG.

**Lemma 22.** *If we choose the stepsize satisfying*

$$\eta \leq \min\left\{\frac{b^2\lambda}{6400L^2}, \frac{1}{8L}\right\}$$

*then the following inequality holds $\forall t \geq 0$:*

$$\begin{aligned}
\mathbb{E}f(w_t) - f(w^*) \leq &\left(1 - \sqrt{\frac{\lambda\eta}{2}}\right)^t \left(V_0(w^*) - f(w^*)\right) \\
&+ \mathbb{E}\left[\sum_{k=1}^{t}\left(1 - \sqrt{\frac{\lambda\eta}{2}}\right)^{t-k}\left(\frac{1}{4\lambda^2\eta^3}\left\|\xi_k\right\|^2 + \frac{8\eta L}{b}(f(w_{k-1}) + f(\widetilde{w}_{s-1}) - 2f(w^*))\right)\right],
\end{aligned}$$

### B.2 Proof of Lemma 16

*Proof.* Since

$$\begin{aligned}
\mathbb{E}\left[f(w_t)\right] &\leq \min_{w}\left\{\mathbb{E}\left[V_t(w)\right]\right\} + \mathbb{E}\left[\delta_t\right] \\
&\leq \mathbb{E}\left[V_t(w^*)\right] + \mathbb{E}\left[\delta_t\right] \\
&\leq (1-\theta_t)f(w^*) + \theta_t V_0(w^*) + \mathbb{E}\left[\delta_t\right],
\end{aligned}$$

subtracting both sides by $f(w^*)$ finishes the proof. □

### B.3 Proof of Lemma 17

*Proof.* We proceed the proof by induction. When $t = 0$, we have

$$V_0(w) = (1-\theta_0)f(w) + \theta_0 V_0(w) = V_0(w),$$

suppose

$$\mathbb{E}\left[V_{t-1}(w)\right] \leq (1-\theta_{t-1})f(w) + \theta_{t-1}V_0(w)$$

is true for some $t \geq 1$, then for $\mathbb{E}[V_t(w)]$, since $\mathbb{E}[v_t] = \nabla f(y_{t-1})$,we have

$$\begin{aligned}
\mathbb{E}[V_t(w)] =&(1-\alpha_{t-1})\mathbb{E}[V_{t-1}(w)] \\
&+ \alpha_{t-1}\mathbb{E}\left[\left(\frac{1}{b}\sum_{i\in B_t} f_i(y_{t-1}) + \left\langle v_t - \frac{\xi_t}{\eta}, w - y_{t-1}\right\rangle + \frac{\lambda}{4}\|w-y_{t-1}\|^2 - \frac{\|\xi_t\|^2}{\lambda\eta^2}\right)\right] \\
=&(1-\alpha_{t-1})\mathbb{E}[V_{t-1}(w)] \\
&+ \alpha_{t-1}\left(f(y_{t-1}) + \langle \nabla f(y_{t-1}), w-y_{t-1}\rangle + \frac{\lambda}{4}\|w-y_{t-1}\|^2 - \left\langle \frac{\xi_t}{\eta}, w-y_{t-1}\right\rangle - \frac{\|\xi_t\|^2}{\lambda\eta^2}\right) \\
\overset{①}{\leq}&(1-\alpha_{t-1})\mathbb{E}[V_{t-1}(w)] + \alpha_{t-1}\left(f(y_{t-1}) + \langle \nabla f(y_{t-1}), w-y_{t-1}\rangle + \frac{\lambda}{2}\|w-y_{t-1}\|^2\right) \\
\overset{②}{\leq}&(1-\alpha_{t-1})(1-\theta_{t-1})f(w) + (1-\alpha_{t-1})\theta_{t-1}V_0(w) + \alpha_{t-1}f(w) \\
=&(1-(1-\alpha_{t-1})\theta_{t-1})f(w) + (1-\alpha_{t-1})\theta_{t-1}V_0(w) \\
=&(1-\theta_t)f(w) + \theta_t V_0(w),
\end{aligned}$$

which concludes the proof, where at step ① we used the inequality $-\langle v_1, v_2\rangle - \|v_2\|^2 \leq \frac{\|v_1\|^2}{4}$ ,at step ② we used the inductive hypothesis and the fact that $f(w)$ is $\lambda$-strongly convex. □

### B.4 Proof of Lemma 18

*Proof.* We proceed the proof by induction, when $t = 0$ this is true by construction. Suppose for some $t \geq 1$ the following holds:

$$V_{t-1}(w) = V_{t-1}^* + \frac{\lambda}{4}\|w - z_{t-1}\|^2,$$

then at time $t$, we have

$$\begin{aligned}
V_t(w) =&(1-\alpha_{t-1})V_{t-1}(w) \\
&+ \alpha_{t-1}\left(\frac{1}{b}\sum_{i\in B_t} f_i(y_{t-1}) + \left\langle v_t - \frac{\xi_t}{\eta}, w - y_{t-1}\right\rangle + \frac{\lambda}{4}\|w-y_{t-1}\|^2 - \frac{\|\xi_t\|^2}{\lambda\eta^2}\right) \\
=&(1-\alpha_{t-1})V_{t-1}^* + \frac{(1-\alpha_{t-1})\lambda}{4}\|w-z_{t-1}\|^2 + \frac{\alpha_{t-1}\lambda}{4}\|w-y_{t-1}\|^2 + \alpha_{t-1}\left\langle v_t - \frac{\xi_t}{\eta}, w\right\rangle \\
&+ \frac{\alpha_{t-1}}{b}\sum_{i\in B_t} f_i(y_{t-1}) - \alpha_{t-1}\left\langle v_t - \frac{\xi_t}{\eta}, y_{t-1}\right\rangle - \frac{\alpha_{t-1}\|\xi_t\|^2}{\lambda\eta^2},
\end{aligned}$$

by first order optimality condition, it is clear that the minimizer of $V_t(w)$: $z_t$ satisfies the following:

$$(1-\alpha_{t-1})\lambda(z_t - z_{t-1}) + \alpha_{t-1}\lambda(z_t - y_{t-1}) + 2\alpha_{t-1}\left(v_t - \frac{\xi_t}{\eta}\right) = 0,$$

from above we obtain the recursive form of $z_t$ as (40) defines. Plug in (40) to (38) we get

$$
\begin{aligned}
V_t^* =& V_t(z_t) = (1-\alpha_{t-1})V_{t-1}^* + \frac{(1-\alpha_{t-1})\lambda}{4}\left\|z_t - z_{t-1}\right\|^2 + \frac{\alpha_{t-1}\lambda}{4}\left\|z_t - y_{t-1}\right\|^2 + \alpha_{t-1}\left\langle v_t - \frac{\xi_t}{\eta}, z_t\right\rangle \\
&+ \frac{\alpha_{t-1}}{b}\sum_{i\in B_t} f_i(y_{t-1}) - \alpha_{t-1}\left\langle v_t - \frac{\xi_t}{\eta}, y_{t-1}\right\rangle - \frac{\alpha_{t-1}\left\|\xi_t\right\|^2}{\lambda\eta^2} \\
=& (1-\alpha_{t-1})V_{t-1}^* + \frac{(1-\alpha_{t-1})\lambda}{4}\left\|\alpha_{t-1}(y_{t-1} - z_{t-1}) - \frac{2\alpha_{t-1}}{\lambda}\left(v_t - \frac{\xi_t}{\eta}\right)\right\|^2 + \frac{\alpha_{t-1}}{b}\sum_{i\in B_t} f_i(y_{t-1}) \\
&+ \frac{\alpha_{t-1}\lambda}{4}\left\|(1-\alpha_{t-1})(z_{t-1} - y_{t-1}) - \frac{2\alpha_{t-1}}{\lambda}\left(v_t - \frac{\xi_t}{\eta}\right)\right\|^2 - \frac{\alpha_{t-1}\left\|\xi_t\right\|^2}{\lambda\eta^2} \\
&+ \alpha_{t-1}\left\langle v_t - \frac{\xi_t}{\eta}, (1-\alpha_{t-1})(z_{t-1} - y_{t-1}) - \frac{2\alpha_{t-1}}{\lambda}\left(v_t - \frac{\xi_t}{\eta}\right)\right\rangle \\
=& (1-\alpha_{t-1})V_{t-1}^* + \frac{\alpha_{t-1}(1-\alpha_{t-1})\lambda}{4}\left\|z_{t-1} - y_{t-1}\right\|^2 + \alpha_{t-1}(1-\alpha_{t-1})\left\langle v_t - \frac{\xi_t}{\eta}, z_{t-1} - y_{t-1}\right\rangle \\
&- \frac{\alpha_{t-1}^2}{\lambda}\left\|v_t\right\|^2 + \frac{\alpha_{t-1}}{b}\sum_{i\in B_t} f_i(y_{t-1}) - \frac{(\alpha_{t-1} + \alpha_{t-1}^2)\left\|\xi_t\right\|^2}{\lambda\eta^2} + \frac{2\alpha_{t-1}^2}{\lambda\eta}\langle \xi_t, v_t\rangle,
\end{aligned}
$$

which verified (41). □

### B.5 Proof of Lemma 19

*Proof.* We prove by induction, when $t = 0$ it is obviously true, suppose it is true for $t-1, \forall t \geq 1$, i.e.

$$
z_{t-1} - y_{t-1} = \frac{1}{\alpha}(y_{t-1} - w_{t-1}),
$$

for iteration $t$, based on (40) we have

$$
\begin{aligned}
z_t - y_t =& (1-\alpha)z_{t-1} + \alpha y_{t-1} - \frac{2\alpha}{\lambda}\left(v_t - \frac{\xi_t}{\eta}\right) - y_t \\
=& (1-\alpha)(z_{t-1} - y_{t-1}) + y_{t-1} - \frac{2\alpha}{\lambda}\left(v_t - \frac{\xi_t}{\eta}\right) - y_t \\
\overset{①}{=}& \frac{1-\alpha}{\alpha}(y_{t-1} - w_{t-1}) + y_{t-1} - \frac{2\alpha}{\lambda}\left(v_t - \frac{\xi_t}{\eta}\right) - y_t \\
=& \frac{1}{\alpha}\left(y_{t-1} - \frac{2\alpha^2}{\lambda}\left(v_t - \frac{\xi_t}{\eta}\right)\right) - \frac{1-\alpha}{\alpha}w_{t-1} - y_t \\
\overset{②}{=}& \frac{1}{\alpha}w_t - \frac{1-\alpha}{\alpha}w_{t-1} - y_t \\
\overset{③}{=}& \frac{1}{\alpha}\left(w_t + \frac{1-\alpha}{1+\alpha}(w_t - w_{t-1}) - w_t\right) \\
\overset{④}{=}& \frac{1}{\alpha}(y_t - w_t),
\end{aligned}
$$

which concludes the proof, where step ① used the inductive hypothesis, step ② used the update rule of $w_t$ in Algorithm 3, step ③ and ④ used the update rule of $y_t$ in Algorithm 3. □

### B.6 Proof of Lemma 20

*Proof.* We prove by induction, when $t = 0$, it is true that

$$f(w_0) \le V_0^* = f(z_0) = f(w_0),$$

suppose

$$\mathbb{E}\left[f(w_{t-1})\right] \le \mathbb{E}\left[V_{t-1}^*\right] + \mathbb{E}\left[\delta_{t-1}\right], \tag{42}$$

for some $t \ge 1$, then based on smoothness, we know

$$f(w_t) \le f(y_{t-1}) + \langle \nabla f(y_{t-1}), w_t - y_{t-1} \rangle + \frac{L}{2}\left\|w_t - y_{t-1}\right\|^2,$$

thus

$$\begin{aligned}
V_t^* =& (1-\alpha_{t-1})V_{t-1}^* + \frac{\alpha_{t-1}(1-\alpha_{t-1})\lambda}{4}\left\|z_{t-1} - y_{t-1}\right\|^2 \\
&+ \alpha_{t-1}(1-\alpha_{t-1})\left\langle v_t - \frac{\xi_t}{\eta}, z_{t-1} - y_{t-1}\right\rangle \\
&- \frac{\alpha_{t-1}^2}{\lambda}\left\|v_t\right\|^2 + \frac{\alpha_{t-1}}{b}\sum_{i\in B_t} f_i(y_{t-1}) - \frac{(\alpha_{t-1} + \alpha_{t-1}^2)\left\|\xi_t\right\|^2}{\lambda\eta^2} + \frac{2\alpha_{t-1}^2}{\lambda\eta}\langle \xi_t, v_t\rangle.
\end{aligned}$$

Hence we have

$$\begin{aligned}
&\mathbb{E}\left[f(w_t) - V_t^*\right] \le \mathbb{E}\left[f(y_{t-1}) + \langle \nabla f(y_{t-1}), w_t - y_{t-1}\rangle + \frac{L}{2}\left\|w_t - y_{t-1}\right\|^2\right] \\
&- \mathbb{E}\left[(1-\alpha)V_{t-1}^* - \frac{\alpha(1-\alpha)\lambda}{4}\left\|z_{t-1} - y_{t-1}\right\|^2 - \alpha(1-\alpha)\left\langle v_t - \frac{\xi_t}{\eta}, z_{t-1} - y_{t-1}\right\rangle\right] \\
&+ \mathbb{E}\left[\frac{\alpha^2}{\lambda}\left\|v_t\right\|^2 - \frac{\alpha}{b}\sum_{i\in B_t} f_i(y_{t-1}) + \frac{(\alpha+\alpha^2)\left\|\xi_t\right\|^2}{\lambda\eta^2} - \frac{2\alpha^2}{\lambda\eta}\langle \xi_t, v_t\rangle\right] \\
&\overset{①}{=} \mathbb{E}\left[(1-\alpha)(f(y_{t-1}) - V_{t-1}^* + \langle v_t, w_{t-1} - y_{t-1}\rangle) - \langle \nabla f(y_{t-1}), \eta v_t - \xi_t\rangle + \frac{\eta}{2}\left\|v_t\right\|^2\right] \\
&+ \mathbb{E}\left[\frac{L}{2}\left\|\eta v_t - \xi_t\right\|^2 - \frac{(1-\alpha)\lambda}{4\alpha}\left\|y_{t-1} - w_{t-1}\right\|^2 + \frac{(\alpha+\alpha^2)\left\|\xi_t\right\|^2}{\lambda\eta^2}\right] \\
&+ \mathbb{E}\left[-\frac{2\alpha^2}{\lambda\eta}\langle \xi_t, v_t\rangle - (1-\alpha)\left\langle \frac{\xi_t}{\eta}, w_{t-1} - y_{t-1}\right\rangle\right] \\
&\overset{②}{\le} \mathbb{E}\left[(1-\alpha)(f(y_{t-1}) - V_{t-1}^* + \langle v_t, w_{t-1} - y_{t-1}\rangle) - \langle \nabla f(y_{t-1}), \eta v_t - \xi_t\rangle + \frac{9\eta}{16}\left\|v_t\right\|^2\right] \\
&+ \mathbb{E}\left[\frac{L}{2}\left\|\eta v_t - \xi_t\right\|^2 + \left(\frac{2\alpha}{\lambda\eta^2} + \frac{16\alpha^4}{\lambda^2\eta^3} + \frac{(1-\alpha)}{2\lambda\eta^2}\right)\left\|\xi_t\right\|^2\right] \\
&+ \mathbb{E}\left[\left(\frac{(1-\alpha)\lambda}{2} - \frac{(1-\alpha)\lambda}{4\alpha}\right)\left\|w_{t-1} - y_{t-1}\right\|^2\right]
\end{aligned} \tag{43}$$

where at step ① we used the fact of updating rule

$$w_t = y_{t-1} - \eta v_t + \xi_t,$$

as well as Lemma 19, and $\mathbb{E}\left[\frac{1}{b}\sum_{i\in B_t} f_i(y_{t-1})\right] = f(y_{t-1})$; and at step ② we used the following two inequalities:

$$\left|(1-\alpha)\left\langle \frac{\xi_t}{\eta}, w_{t-1}-y_{t-1}\right\rangle\right| \le \frac{(1-\alpha)\lambda}{2}\left\|y_{t-1}-w_{t-1}\right\|^2 + \frac{(1-\alpha)}{2\lambda\eta^2}\left\|\xi_t\right\|^2,$$
$$\left|\frac{2\alpha^2}{\lambda\eta}\langle \xi_t, v_t\rangle\right| \le \frac{\eta}{16}\left\|v_t\right\|^2 + \frac{16\alpha^4}{\lambda^2\eta^3}\left\|\xi_t\right\|^2$$

Since

$$\begin{aligned}\mathbb{E}\left[f(y_{t-1}) - V_{t-1}^* + \langle v_t, w_{t-1}-y_{t-1}\rangle\right] &= \mathbb{E}\left[f(y_{t-1}) - V_{t-1}^* + \langle \nabla f(y_{t-1}), w_{t-1}-y_{t-1}\rangle\right]\\ &\overset{①}{\le} \mathbb{E}\left[f(w_{t-1}) - \frac{\lambda}{2}\left\|w_{t-1}-y_{t-1}\right\|^2 - V_{t-1}^*\right]\\ &\overset{②}{\le} \mathbb{E}\left[\delta_{t-1} - \frac{\lambda}{2}\left\|w_{t-1}-y_{t-1}\right\|^2\right], \end{aligned} \tag{44}$$

where at step ① we used the $\lambda$-strong convexity of $f(w)$, at step ② we used the inductive hypothesis (42). Combining (43) and (44) we obtain

$$\begin{aligned}&\mathbb{E}\left[f(w_t) - V_t^*\right] \le (1-\alpha)\delta_{t-1} - \mathbb{E}\left[\frac{(1-\alpha)\lambda}{4\alpha}\left\|y_{t-1}-w_{t-1}\right\|^2\right]\\ &+\mathbb{E}\left[\frac{L}{2}\left\|\eta v_t - \xi_t\right\|^2 + \frac{9\eta}{16}\left\|v_t\right\|^2 - \langle \nabla f(y_{t-1}), \eta v_t - \xi_t\rangle + \left(\frac{2\alpha}{\lambda\eta^2} + \frac{16\alpha^4}{\lambda^2\eta^3} + \frac{(1-\alpha)}{2\lambda\eta^2}\right)\left\|\xi_t\right\|^2\right]. \end{aligned} \tag{45}$$

Next we bound the third term on the right hand side, since

$$\begin{aligned}&\frac{L}{2}\left\|\eta v_t - \xi_t\right\|^2 + \frac{9\eta}{16}\left\|v_t\right\|^2 - \langle \nabla f(y_{t-1}), \eta v_t + \xi_t\rangle\\ &\overset{①}{\le} L\eta^2\left\|v_t\right\|^2 + L\left\|\xi_t\right\|^2 + \frac{9\eta}{16}\left\|v_t\right\|^2 - \langle \nabla f(y_{t-1}), \eta v_t - \xi_t\rangle\\ &= \frac{\eta}{2}\left\|v_t - \nabla f(y_{t-1})\right\|^2 - \frac{\eta}{2}\left\|\nabla f(y_{t-1})\right\|^2\\ &+\left(L\eta^2 + \frac{\eta}{16}\right)\left\|v_t\right\|^2 + \langle \nabla f(y_{t-1}), \xi_t\rangle + L\left\|\xi_t\right\|^2\\ &\overset{②}{\le} \frac{\eta}{2}\left\|v_t - \nabla f(y_{t-1})\right\|^2 - \frac{\eta}{2}\left\|\nabla f(y_{t-1})\right\|^2\\ &+\left(2L\eta^2 + \frac{\eta}{8}\right)\left(\left\|v_t - \nabla f(y_{t-1})\right\|^2 + \left\|\nabla f(y_{t-1})\right\|^2\right)\\ &+\frac{\eta}{8}\left\|\nabla f(y_{t-1})\right\|^2 + \frac{2}{\eta}\left\|\xi_t\right\|^2 + L\left\|\xi_t\right\|^2, \end{aligned}$$

where at step ① we used $\|v_1+v_2\|^2 \le 2(\|v_1\|^2 + \|v_2\|^2)$, and at step ② we used it again and also the inequality $2|\langle v_1, v_2\rangle| \le \|v_1\|^2 + \|v_2\|^2$. Thus we know when $\eta \le \frac{1}{8L}$, we have

$$\frac{L}{2}\left\|\eta v_t - \xi_t\right\|^2 + \frac{9\eta}{16}\left\|v_t\right\|^2 - \langle \nabla f(y_{t-1}), \eta v_t + \xi_t\rangle \tag{46}$$
$$\begin{aligned}&\le \left(\frac{\eta}{2} + 2L\eta^2 + \frac{\eta}{8}\right)\left\|v_t - \nabla f(y_{t-1})\right\|^2 + \left(L + \frac{2}{\eta}\right)\left\|\xi_t\right\|^2\\ &-\left(\frac{\eta}{2} - \frac{\eta}{8} - \frac{\eta}{8} - 2L\eta^2\right)\left\|\nabla f(y_{t-1})\right\|^2\\ &\le \eta\left\|v_t - \nabla f(y_{t-1})\right\|^2 + \frac{3}{\eta}\left\|\xi_t\right\|^2. \end{aligned} \tag{47}$$

Combining (45) and (47) we obtain

$$\begin{aligned}\mathbb{E}\left[f(w_t)-V_t^*\right] \leq & (1-\alpha)\delta_{t-1}+\mathbb{E}\left[\eta\left\|v_t-\nabla f(y_{t-1})\right\|^2-\frac{(1-\alpha)\lambda}{4\alpha}\left\|y_{t-1}-w_{t-1}\right\|^2\right] \\ & +\mathbb{E}\left[\left(\frac{2\alpha}{\lambda\eta^2}+\frac{16\alpha^4}{\lambda^2\eta^3}+\frac{(1-\alpha)}{2\lambda\eta^2}+\frac{3}{\eta}\right)\left\|\xi_t\right\|^2\right],\end{aligned} \tag{48}$$

also since

$$\alpha=\sqrt{\frac{\lambda\eta}{2}} \leq \frac{1}{4}\sqrt{\frac{\lambda}{L}} \leq \frac{1}{4},$$

we know

$$\frac{2\alpha}{\lambda\eta^2}+\frac{16\alpha^4}{\lambda^2\eta^3}+\frac{(1-\alpha)}{2\lambda\eta^2}+\frac{3}{\eta} \leq \frac{1}{16\lambda^2\eta^3}+\frac{1}{16\lambda^2\eta^3}+\frac{1}{16\lambda^2\eta^3}+\frac{3}{64\lambda^2\eta^3} \leq \frac{1}{4\lambda^2\eta^3},$$

combining the inequality above and (48), then substituting $\alpha=\sqrt{\frac{\lambda\eta}{2}}$ finishes the proof. □

## B.7 Proof of Lemma 21

*Proof.* The proof follows the strategy in (Johnson and Zhang, 2013; Nitanda, 2014). First, based on the minibatch sampling (here for simplicity we only consider sampling with replacement, for sampling without replacement, the bound below can be tightened by a factor of $\frac{n-b}{n-1}$, see, e.g. Section 2.8 of (Lohr, 2009)) we know

$$\begin{aligned}\mathbb{E}\left[\left\|v_t-\nabla f(y_{t-1})\right\|^2\right] = & \frac{1}{b^2}\mathbb{E}\left[\sum_{j\in B_t}\left\|\nabla f_j(y_{t-1})-\nabla f(y_{t-1})-(\nabla f_j(\widetilde{w}_{s-1})-\nabla f(\widetilde{w}_{s-1}))\right\|^2\right] \\ \overset{①}{\leq} & \frac{1}{b^2}\mathbb{E}\left[\sum_{j\in B_t}\left\|\nabla f_j(y_{t-1})-\nabla f_j(\widetilde{w}_{s-1})\right\|^2\right] \\ \overset{②}{\leq} & \frac{2}{b^2}\mathbb{E}\left[\sum_{j\in B_t}\left(\left\|\nabla f_j(y_{t-1})-\nabla f_j(w_{t-1})\right\|^2+\left\|\nabla f_j(\widetilde{w}_{s-1})-\nabla f_j(w_{t-1})\right\|^2\right)\right] \\ \overset{③}{\leq} & \frac{4}{b^2}\mathbb{E}\left[\sum_{j\in B_t}\left(\left\|\nabla f_j(\widetilde{w}_{s-1})-\nabla f_j(w^*)\right\|^2\right)+\left\|\nabla f_j(w^*)-\nabla f_j(w_{t-1})\right\|^2\right)\right] \\ & +\frac{2L^2}{b}\left\|y_{t-1}-w_{t-1}\right\|^2 \\ \overset{④}{\leq} & \frac{8L}{b}\left(f(w_{t-1})-f(w^*)+f(\widetilde{w}_{s-1})-f(w^*)\right)+\frac{2L^2}{b}\left\|y_{t-1}-w_{t-1}\right\|^2,\end{aligned}$$

where at step ① we used $\mathbb{E}\left\|v-\mathbb{E}v\right\|^2 \leq \mathbb{E}\left\|v\right\|^2$, at step ② we used $\left\|v_1+v_2\right\|^2 \leq 2(\left\|v_1\right\|^2+\left\|v_2\right\|^2)$, at step ③ we used it again along with the $L$-smoothness of $f_j(w)$, at step ④ we used standard results in SVRG analysis, e.g. Lemma 1 of (Xiao and Zhang, 2014). Substituting above into the term $\mathbb{E}\left[\eta\left\|v_t-\nabla f(y_{t-1})\right\|^2-\frac{(1-\sqrt{\lambda\eta/2})\lambda}{\sqrt{2\lambda\eta}}\left\|y_{t-1}-w_{t-1}\right\|^2\right]$ concludes the proof. □

## B.8 Proof of Lemma 22

*Proof.* Combining Lemma 20 and 21, we have

$$\begin{aligned}
\delta_t &= \left(1-\sqrt{\frac{\lambda\eta}{2}}\right)\delta_{t-1} \\
&\quad + \left(\eta\left\|v_t-\nabla f(y_{t-1})\right\|^2 - \frac{(1-\sqrt{\lambda\eta/2})\lambda}{\sqrt{2\lambda\eta}}\left\|y_{t-1}-w_{t-1}\right\|^2 + \frac{1}{4\lambda^2\eta^3}\left\|\xi_t\right\|^2\right) \qquad (49)\\
&\le \left(1-\sqrt{\frac{\lambda\eta}{2}}\right)\delta_{t-1} + \left(\frac{1}{4\lambda^2\eta^3}\left\|\xi_t\right\|^2 + \frac{8\eta L}{b}(f(w_{t-1})+f(\widetilde{w}_{s-1})-2f(w^*))\right) \\
&\quad + \left(\frac{2\eta L^2}{b} - \frac{(1-\sqrt{\lambda\eta/2})\lambda}{\sqrt{2\lambda\eta}}\right)\left\|y_{t-1}-w_{t-1}\right\|^2. \qquad (50)
\end{aligned}$$

First we verify the factor before $\|y_{t-1}-w_{t-1}\|^2$ is non positive, since $1-\eta\lambda \ge 1-\frac{\lambda}{8L} \ge \frac{1}{2}$, then

$$\frac{2\eta L^2}{b}\cdot\frac{\sqrt{2\lambda\eta}}{(1-\sqrt{\eta\lambda/2})\lambda} \le \frac{4\eta^{3/2}L^2}{b\sqrt{\lambda}} \le \min\left\{\frac{b^2\lambda}{128000L}, \frac{L^{1/2}}{4b\sqrt{\lambda}}\right\} \le 1,$$

where the last inequality is true because

$$\frac{b^2\lambda}{128000L}\left(\frac{L^{1/2}}{4b\sqrt{\lambda}}\right)^2 \le 1,$$

thus we must have $\min\left\{\frac{b^2\lambda}{128000L}, \frac{L^{1/2}}{4b\sqrt{\lambda}}\right\} \le 1$ otherwise it leads to a contradiction. Thus $\frac{2\eta L^2}{b} - \frac{(1-\sqrt{\eta\lambda/2})\lambda}{\sqrt{2\eta\lambda}} \le 0$, combining with (50), we have

$$\delta_t \le \left(1-\sqrt{\frac{\lambda\eta}{2}}\right)\delta_{t-1} + \left(\frac{1}{4\lambda^2\eta^3}\left\|\xi_t\right\|^2 + \frac{8\eta L}{b}(f(w_{t-1})+f(\widetilde{w}_{s-1})-2f(w^*))\right),$$

applying above inequality recursively, we get

$$\delta_t \le \sum_{k=1}^{t}\left(1-\sqrt{\frac{\lambda\eta}{2}}\right)^{t-k}\left(\frac{1}{4\lambda^2\eta^3}\left\|\xi_k\right\|^2 + \frac{8\eta L}{b}(f(w_{k-1})+f(\widetilde{w}_{s-1})-2f(w^*))\right),$$

then applying Lemma 16 and Lemma 20 concludes the proof. □

## B.9 Proof of Theorem 13

*Proof.* We prove this theorem by induction, note that when $b \le 20\sqrt{\frac{2L}{\lambda}}$, it is clear that

$$\frac{8\eta L}{b} \le \frac{b^2\lambda}{6400L^2}\cdot\frac{8L}{b} = \frac{\lambda b}{800L} = \frac{\sqrt{\eta\lambda}}{10},$$

and when $b \ge 20\sqrt{\frac{2L}{\lambda}}$,

$$\frac{8\eta L}{b} = \frac{1}{b} \le \frac{\sqrt{\lambda}}{20\sqrt{2L}} = \frac{\sqrt{\eta\lambda}}{10},$$

thus in both cases we have $\frac{8\eta L}{b} \le \frac{\sqrt{\eta\lambda}}{10}$, thus

$$\begin{aligned}\sum_{k=1}^{t}\left(1-\sqrt{\frac{\lambda\eta}{2}}\right)^{t-k}\left(\frac{8\eta L}{b}\right) \le & \frac{1}{10}\sum_{k=1}^{t}\left(1-\sqrt{\frac{\lambda\eta}{2}}\right)^{t-k}\sqrt{\lambda\eta} \\ \le & \frac{1}{10}\sum_{k=1}^{\infty}\left(1-\sqrt{\frac{\lambda\eta}{2}}\right)^{k}\sqrt{\lambda\eta} \le \frac{\sqrt{2}}{10}.\end{aligned}$$

By Lemma 22 we know

$$\mathbb{E}f(w_t) - f(w^*) \le \left(1-\sqrt{\frac{\lambda\eta}{2}}\right)^{t}(V_0(w^*) - f(w^*)) + \frac{\sqrt{2}(f(\widetilde{w}_{s-1}) - f(w^*))}{10} \tag{51}$$

$$+\mathbb{E}\left[\sum_{k=1}^{t}\left(1-\sqrt{\frac{\lambda\eta}{2}}\right)^{t-k}\left(\frac{1}{4\lambda^2\eta^3}\|\xi_k\|^2 + \frac{\sqrt{\lambda\eta}}{10}(f(w_{k-1}) - f(w^*))\right)\right], \tag{52}$$

also when $\forall k \in [t]$,

$$\|\xi_k\|^2 \le \frac{2\lambda^2\eta^3\sqrt{\lambda\eta}}{15}(f(\widetilde{w}_{s-1}) - f(w^*)), \tag{53}$$

we know

$$\begin{aligned}\sum_{k=1}^{t}\left(1-\sqrt{\frac{\lambda\eta}{2}}\right)^{t-k}\left(\frac{1}{4\lambda^2\eta^3}\|\xi_k\|^2\right) \le & \frac{1}{30}\sum_{k=1}^{t}\left(1-\sqrt{\frac{\lambda\eta}{2}}\right)^{t-k}\sqrt{\lambda\eta}(f(\widetilde{w}_{s-1}) - f(w^*)) \\ \le & \frac{1}{30}\sum_{k=1}^{\infty}\left(1-\sqrt{\frac{\lambda\eta}{2}}\right)^{k}\sqrt{\lambda\eta}(f(\widetilde{w}_{s-1}) - f(w^*)) \\ \le & \frac{(f(\widetilde{w}_{s-1}) - f(w^*))}{20}.\end{aligned} \tag{54}$$

Combining (52) and (54), we have

$$\begin{aligned}\mathbb{E}f(w_t) - f(w^*) \le & \left(1-\sqrt{\frac{\lambda\eta}{2}}\right)^{t}(V_0(w^*) - f(w^*)) + \frac{(f(\widetilde{w}_{s-1}) - f(w^*))}{5} \\ & +\mathbb{E}\left[\sum_{k=1}^{t}\left(1-\sqrt{\frac{\lambda\eta}{2}}\right)^{t-k}\left(\frac{\sqrt{\lambda\eta}}{10}(f(w_{k-1}) - f(w^*))\right)\right],\end{aligned}$$

denote $A_t$ be the right hand side on above inequality, i.e.

$$\begin{aligned}
A_t &= \left(1-\sqrt{\frac{\lambda\eta}{2}}\right)^t (V_0(w^*) - f(w^*)) + \frac{(f(\widetilde{w}_{s-1}) - f(w^*))}{5} \\
&\quad + \mathbb{E}\left[\sum_{k=1}^{t}\left(1-\sqrt{\frac{\lambda\eta}{2}}\right)^{t-k}\left(\frac{\sqrt{\lambda\eta}}{10}(f(w_{k-1}) - f(w^*))\right)\right] \\
&= \left(1-\sqrt{\frac{\lambda\eta}{2}}\right)\left(\left(1-\sqrt{\frac{\lambda\eta}{2}}\right)^{t-1}(V_0(w^*) - f(w^*)) + \frac{(f(\widetilde{w}_{s-1}) - f(w^*))}{5}\right) \\
&\quad + \left(1-\sqrt{\frac{\lambda\eta}{2}}\right)\mathbb{E}\left[\sum_{k=1}^{t-1}\left(1-\sqrt{\frac{\lambda\eta}{2}}\right)^{t-1-k}\left(\frac{\sqrt{\lambda\eta}}{10}(f(w_{k-1}) - f(w^*))\right)\right] \\
&\quad + \frac{\sqrt{\lambda\eta}}{5\sqrt{2}}(f(\widetilde{w}_{s-1}) - f(w^*)) + \frac{\sqrt{\lambda\eta}}{10}(f(w_{t-1}) - f(w^*)) \\
&\overset{①}{\leq} \left(1-\sqrt{\frac{\lambda\eta}{2}}\right)A_{t-1} + \frac{\sqrt{\lambda\eta}}{5\sqrt{2}}A_0 + \frac{\sqrt{\lambda\eta}}{10}A_{t-1} \\
&= \left(1-\frac{9\sqrt{\lambda\eta}}{10}\right)A_{t-1} + \frac{\sqrt{\lambda\eta}}{5\sqrt{2}}A_0,
\end{aligned}$$

where at step ① we used the fact that $f(\widetilde{w}_{s-1}) - f(w^*) \leq A_0$, and $f(w_{t-1}) - f(w^*) \leq A_{t-1}$, applying above inequality recursively we get

$$\begin{aligned}
A_t &\leq \left(1-\frac{9\sqrt{\lambda\eta}}{10}\right)^t A_0 + \frac{\sqrt{\lambda\eta}}{5\sqrt{2}}\sum_{k=0}^{t-1}\left(1-\frac{9\sqrt{\lambda\eta}}{10}\right)^k A_0 \\
&\leq \left(1-\frac{9\sqrt{\lambda\eta}}{10}\right)^t A_0 + \frac{\sqrt{\lambda\eta}}{5\sqrt{2}}\sum_{k=0}^{\infty}\left(1-\frac{9\sqrt{\lambda\eta}}{10}\right)^k A_0 \leq \left(\left(1-\frac{9\sqrt{\lambda\eta}}{10}\right)^t + \frac{2}{9}\right)A_0,
\end{aligned}$$

thus we know when $t \geq \frac{10}{9\sqrt{\lambda\eta}}\log(36)$, we have

$$A_t \leq \left(\frac{1}{36} + \frac{2}{9}\right)A_0 = \frac{1}{4}A_0,$$

also since

$$A_0 = V_0(w^*) - f(w^*) = f(\widetilde{w}_{s-1}) - f(w^*) + \frac{\lambda}{2}\|\widetilde{w}_{s-1} - w^*\| \overset{①}{\leq} 2(f(\widetilde{w}_{s-1}) - f(w^*)),$$

where at step ① we used the $\lambda$-strong convexity of $f(w)$. Combine above two inequality we get when $t \geq \frac{10}{9\sqrt{\lambda\eta}}\log(36)$, the expected objective suboptimality is halved:

$$\mathbb{E}f(w_t) - f(w^*) \leq A_t \leq \frac{1}{4}A_0 \leq \frac{1}{2}(f(\widetilde{w}_{s-1}) - f(w^*)),$$

which concludes the proof. □

## C Collections of Tools in the Analysis

**Lemma 23.** *(**Matrix Bernstein**, Theorem 1.4 of (Tropp, 2012) rephrased) Let $X_1, ..., X_k$ be some independent, self-adjoint random matrices with dimension $d$, and assume each random matrix satisfies:*

$$\mathbb{E}X_k = 0 \quad \text{and} \quad \|X_k\| \leq R \quad \textit{almost surely.}$$

*Then, for all $t \geq 0$,*

$$P\left(\left\|\sum_k X_k\right\| \geq t\right) \leq d \exp\left(\frac{-t^2/2}{\sigma^2 + Rt/3}\right)$$

*where*

$$\sigma^2 := \left\|\sum_k \mathbb{E}(X_k^2)\right\|.$$

**Lemma 24.** *(**Iteration complexity of SVRG**, Corollary 1 of (Xiao and Zhang, 2014) rephrased) If we apply SVRG to any finite-sum optimization objective $f(w)$ where each individual function is $\lambda$-strongly convex and $L$-smooth, then there are universal constant $C$ such that for any target accuracy $\varepsilon$, and success probability $1-\delta$, SVRG is able to find a solution that satisfies $\varepsilon$ objective suboptimality using*

$$C \cdot \left(n + \frac{L}{\lambda}\right) \log\left(\frac{\varepsilon_0}{\delta\varepsilon}\right)$$

*first order oracle calls of individual functions, where $n$ is the total number of individual functions in $f(w)$, $\varepsilon_0$ is the initial objective suboptimality.*

**Lemma 25.** *(**Iteration complexity of catalyst acceleration**, Theorem 3.1 of (Lin et al., 2015a) rephrased) For any $\lambda$-strongly convex and $L$-smooth function $f(w)$, If the minimization step of (20) satisfies*

$$f(w_r) - \frac{\gamma}{2}\|w_r - z_{r-1}\| - \min_w \left(f(w) + \frac{\gamma}{2}\|w - z_{r-1}\|\right) \leq \frac{2}{9}(f(w_0) - f(w^*))\left(1 - \frac{9}{10}\sqrt{\frac{\lambda}{\lambda+\gamma}}\right)^r,$$

*then if we initialize $\nu_0 = \sqrt{\frac{\lambda}{\lambda+\gamma}}$ and set $\nu_r$ such that $\nu_r^2 = (1-\nu_r)\nu_{r-1}^2 + (\lambda\nu_r)/(\lambda+\gamma)$, then the sequences $\{w_r\}$ in Algorithm 2 satisfies*

$$f(w_r) - f(w^*) \leq \frac{800(\lambda+\gamma)}{\lambda}\left(1 - \frac{9}{10}\sqrt{\frac{\lambda}{\lambda+\gamma}}\right)^{r+1}(f(w_0) - f(w^*)).$$